\documentclass[11pt]{article}
\usepackage{graphicx}
\usepackage{bm}
\usepackage{amsmath,amsthm,amssymb,url}
\usepackage{multirow,bigdelim}
\usepackage{amscd}

\usepackage{tikz}
\usetikzlibrary{arrows}

\usepackage{hyperref}

\author{James Cruickshank\thanks{School of Mathematical and Statistical Sciences, University of Galway, Ireland }, Csaba Kir\'aly\thanks{
ELTE E\"otv\"os Lor\'and University and the HUN-REN--ELTE Egerv\'ary Research Group on Combinatorial Optimization, 
P\'azm\'any
P\'eter s\'et\'any 1/C
Budapest, Hungary H-1117},
Bernd Schulze\thanks{School of Mathematical Sciences, Lancaster University, Lancaster, UK},
Shin-ichi Tanigawa\thanks{Department of Business Economics, School of Management, Tokyo University of Science, Japan}
}
\title{On unique realisations of graphs under symmetry constraints}
\date{\today}

\DeclareMathOperator{\rank}{rank}

\renewcommand{\tilde}{\widetilde}

\newtheorem{thm}{Theorem}[section]
\newtheorem{theorem}[thm]{Theorem}
\newtheorem{proposition}[thm]{Proposition}
\newtheorem{lemma}[thm]{Lemma}

\newtheorem{claim}[thm]{Claim}

\newtheorem{conjecture}{Conjecture}
\newtheorem{remark}[thm]{Remark}

\begin{document}
\maketitle

\begin{abstract}
\noindent Global rigidity asks whether a realisation of a graph is uniquely
determined, up to Euclidean congruence, by its edge lengths.
In this paper, we study global rigidity of bar-joint frameworks under
finite point-group symmetry constraints.
A framework is called forced-symmetric globally rigid in $d$-space
if every realisation of the graph in $\mathbb R^d$ with the same
edge lengths and the same prescribed symmetry is congruent to it.
We establish symmetric extensions of two fundamental results on
global rigidity.

Our first result extends Hendrickson's necessary condition to the
forced-symmetric setting. For symmetry-generic frameworks with an
arbitrary finite point group acting freely on the vertices and with
at least $d+2$ vertex orbits, we prove that forced-symmetric global
rigidity implies forced-symmetric redundant rigidity.

Our second result extends the classical characterisation of generic
global rigidity for body-bar frameworks due to Connelly, Jord\'an
and Whiteley. For the symmetric body-bar model considered here,
we show that forced-symmetric global rigidity is a generic property
and obtain a combinatorial characterisation in every dimension.
\end{abstract}


\section{Introduction}

A bar-joint framework is a realisation of a graph in Euclidean space,
with vertices representing joints and edges representing rigid bars.
Local rigidity asks whether the framework is isolated, up to Euclidean
motions, among nearby realisations with the same edge lengths, whereas
global rigidity asks whether it is the unique such realisation among
all realisations. These questions arise in structural engineering,
robotics, and materials science, as well as in sensor network
localisation \cite{aspnes} and the analysis of biomolecular
structures \cite{hermans}.

Symmetry is a prominent feature of many of these structures, from
engineered constructions to biomolecular assemblies. Its effects on
rigidity and flexibility have therefore been studied for many decades;
see \cite{HoDCG} for a survey. Interest in this subject grew rapidly
in the 2010s, as it became clear that gain graphs could be used to
extend many classical rigidity theorems to symmetric
frameworks \cite{bs10,mt13,rossthesis}. By encoding a
symmetric graph through its quotient graph together with group labels
on the edges, this approach makes it possible to express rigidity
conditions in terms of both the graph structure and the prescribed
symmetry. It has led to combinatorial characterisations of local
rigidity for several important classes of symmetric frameworks.

In this paper, we treat symmetry as an additional constraint on the realisation
space: alternative realisations must satisfy both the prescribed edge
lengths and the prescribed symmetry. This is the
\emph{forced-symmetric} setting, which has been studied extensively
for {\em infinitesimal} and {\em local} rigidity
\cite{jkt,mt15,BSsiam,T2,BSWWorbit,bernstein,dewaretal,HoDCG}.
Here, the {\em global} question is
whether these constraints determine a framework uniquely up to
Euclidean congruence. Symmetry can change the global rigidity question substantially.
For example, a framework with a disconnected underlying graph may
be globally rigid among realisations satisfying the prescribed
symmetry, since the symmetry can constrain the relative positions
of its connected components.

The theory of global rigidity in this setting is much less developed
than that of local rigidity. Recent progress for infinite periodic
frameworks \cite{kst,kks,dstt} has shown that global rigidity can be
analysed effectively when periodicity is imposed as a constraint. In
this paper, we develop the theory for finite frameworks with
point-group symmetry, extending two fundamental results from
classical global rigidity theory to this setting.

Passing from periodicity to point-group symmetry introduces further
difficulties. In the full-rank periodic setting with a fixed lattice
representation, the Euclidean isometries commuting with the lattice
translations are simply translations. For point-group symmetry,
the corresponding group of isometries is the centraliser of the
point group in the Euclidean group, and its structure depends on
the prescribed representation. Consequently, the geometry of the
realisation space modulo symmetry-preserving congruences varies
with the symmetry. The compactness and quotient arguments underlying
classical global rigidity proofs therefore require additional care
in this setting.

Our first main result is a forced-symmetric analogue of Hendrickson's
necessary condition for global rigidity \cite{hend}. In the classical
setting, a generic globally rigid framework with sufficiently many
vertices must be redundantly rigid, that is, remain locally rigid
after the deletion of any edge. This provides a fundamental link
between global and local rigidity. For generic symmetric frameworks
under suitable assumptions, we prove that
forced-symmetric global rigidity implies forced-symmetric redundant
rigidity. This result allows the existing theory
of forced-symmetric local rigidity to provide necessary conditions
for global rigidity. Our proof first establishes an abstract version
of Hendrickson's argument under a compactness assumption, and then
shows that global rigidity supplies the boundedness needed to verify
this assumption in the symmetric setting.

A natural next question is when this necessary condition is also
sufficient. Hendrickson conjectured that redundant rigidity, together
with $(d+1)$-vertex-connectivity, characterises generic global rigidity
in $\mathbb{R}^d$. Connelly subsequently disproved this conjecture
in dimensions $d\geq 3$ \cite{Con91}, while the planar case
was later established by Jackson and Jord\'an \cite{jj}.
Identifying substantial classes of frameworks for which
Hendrickson's necessary conditions are sufficient has since
remained an important direction of research. A notable
example is the theorem of Connelly, Jord\'an and Whiteley
\cite{cjw}, which states that, in every dimension, a generic body-bar framework is globally
rigid if and only if it is redundantly rigid.

Our second main result establishes a forced-symmetric counterpart of
this theorem. Body-bar frameworks consist of full-dimensional rigid
bodies connected by bars and provide a standard model in engineering. 
For the symmetric body-bar model considered here, we characterise
generic global rigidity by forced-symmetric
redundant rigidity, connectivity of the quotient graph, and
equality of the centralisers associated with a covering component
and with the whole group. When the covering graph is connected,
the centraliser condition is automatic,
and the global rigidity can be characterised purely by redundant rigidity. 
Combining our result with the existing combinatorial characterisation of
forced-symmetric local rigidity for body-bar frameworks \cite{T2}
yields a polynomial-time checkable characterisation of
forced-symmetric global rigidity in arbitrary dimension.

The paper is organised as follows. Section~\ref{sec:prelim}
recalls the necessary background on global
rigidity and symmetric frameworks.
Section~\ref{sec:symgr} introduces forced-symmetric local and global rigidity
and their formulation in terms of quotient gain graphs, and gives a rank
characterisation of generic local rigidity.
Section~\ref{sec:Hen} establishes the forced-symmetric analogue of
Hendrickson's necessary condition, and Section~\ref{sec:bb} proves
the characterisation for body-bar frameworks.
Section~\ref{sec:future} discusses further questions.

\section[Preliminaries]{Preliminaries}\label{sec:prelim}

In this section, we provide basic definitions of global rigidity.
For $x,y\in\mathbb R^d$, we write
$\langle x,y\rangle=x^\top y$,
$\|x\|^2=x^\top x=\sum_{k=1}^d x_k^2$.

\subsection[Background on global rigidity]{Background on global rigidity}\label{sec:backgr}

A \emph{bar-joint framework} (or simply a \emph{framework}) in
$\mathbb R^d$ is a pair $(\tilde G,\tilde p)$, where
$\tilde G=(\tilde V,\tilde E)$ is a simple graph and
$\tilde p:\tilde V\to\mathbb R^d$ is a map.
We also call $(\tilde G,\tilde p)$ a $d$-dimensional
\emph{realisation} of the \emph{underlying graph} $\tilde G$,
and write $\tilde p_i=\tilde p(i)$.
For $\tilde V=\{1,\ldots,n\}$, we identify $\tilde p$ with a vector
in $\mathbb R^{dn}$, called a \emph{configuration}.

For a fixed ordering of the edges of $\tilde G$, we define the
\emph{edge function} $f_{\tilde G}:\mathbb R^{dn}\to
\mathbb R^{|\tilde E|}$ by
\[
f_{\tilde G}(\tilde p_1,\ldots,\tilde p_n)
=\bigl(\ldots,\|\tilde p_i-\tilde p_j\|^2,\ldots\bigr),
\]
where $ij\in\tilde E$. 
Its fibre $f_{\tilde G}^{-1}(f_{\tilde G}(\tilde p))$ consists of
all configurations in $\mathbb R^{dn}$ with the same squared edge
lengths as $\tilde p$. 
Let  $K_n$ be the complete graph on $\tilde V$.
Then $f_{K_n}^{-1}(f_{K_n}(\tilde p))
\subseteq f_{\tilde G}^{-1}(f_{\tilde G}(\tilde p))$.

Two frameworks $(\tilde G,\tilde p)$ and $(\tilde G,\tilde q)$ are
\emph{equivalent}, written
$(\tilde G,\tilde p)\simeq(\tilde G,\tilde q)$, if
$f_{\tilde G}(\tilde p)=f_{\tilde G}(\tilde q)$.
They are \emph{congruent}, written
$(\tilde G,\tilde p)\equiv(\tilde G,\tilde q)$, if
$f_{K_n}(\tilde p)=f_{K_n}(\tilde q)$.

We write $O(\mathbb R^d)$ for the orthogonal group of $\mathbb R^d$.
An \emph{isometry} of $\mathbb R^d$ is a map $x\mapsto Ax+t$,
where $A\in O(\mathbb R^d)$ and $t\in\mathbb R^d$.
We denote the group of isometries of $\mathbb R^d$ by
$\mathrm{Euc}(d)$.

Two frameworks $(\tilde G,\tilde p)$ and $(\tilde G,\tilde q)$
are congruent if and only if there exists an isometry
$h\in\mathrm{Euc}(d)$ such that $\tilde q_i=h(\tilde p_i)$
for every $i\in\tilde V$.

A framework $(\tilde G,\tilde p)$ is \emph{(locally) rigid} in
$\mathbb R^d$ if there exists an open neighbourhood $U$ of
$\tilde p$ in $\mathbb R^{dn}$ such that, for every $\tilde q\in U$,
equivalence to $(\tilde G,\tilde p)$ implies congruence to it.
Here and below, neighbourhoods are taken in the usual topology of
$\mathbb R^{dn}$.
The framework $(\tilde G,\tilde p)$ is \emph{globally rigid} in
$\mathbb R^d$ if, for every $\tilde q\in\mathbb R^{dn}$,
$(\tilde G,\tilde p)\simeq(\tilde G,\tilde q)$
implies $(\tilde G,\tilde p)\equiv(\tilde G,\tilde q)$.

A configuration $\tilde p\in\mathbb R^{dn}$ is \emph{generic} if
its coordinates are algebraically independent over $\mathbb Q$.
A framework $(\tilde G,\tilde p)$ is \emph{generic} if its
configuration is generic.

\subsection[Symmetric graphs and frameworks]{Symmetric graphs and frameworks}\label{sec:symbasic}

Let $\tilde G$ be a finite simple graph. An {\em action} of a finite group $\Gamma$ on $\tilde G$ is a group homomorphism $\theta:\Gamma \rightarrow {\rm Aut}(\tilde G)$, where ${\rm Aut}(\tilde G)$ denotes the automorphism group of $\tilde G$. We say that  $\tilde G$ is {\em $\Gamma$-symmetric} (with respect to $\theta$) if $\Gamma$ acts on $\tilde G$ by $\theta$.  In the following we will frequently
omit specifying the action $\theta$ if it is clear from the context.
We then denote $\theta(\gamma)(i)$ by  $\gamma i$. For simplicity we will assume that $\Gamma$  acts freely on $\tilde V$ via $\theta$, i.e. $\theta(\gamma)(i)\neq i$ for every $\gamma\in\Gamma\setminus\{1\}$ and every $i\in\tilde V$.

For a $\Gamma$-symmetric graph $\tilde G=(\tilde V,\tilde E)$, the orbit of a vertex $i$ (respectively an edge $e$) is denoted by $\Gamma i$ (respectively $\Gamma e$). Thus, $\Gamma i=\{\gamma i|\, \gamma\in \Gamma\}$ and $\Gamma e=\{\gamma e|\, \gamma\in \Gamma\}$. The collection of all vertex
orbits (respectively, edge orbits) is denoted by $V$ (respectively $E$). The \emph{quotient graph} $\tilde G /\Gamma$ of $\tilde G$ is a multigraph with vertex set $V$, edge set $E$ and  incidence relation satisfying $\Gamma e=\Gamma i \Gamma j$ if some (equivalently every) edge in $\Gamma e$ is incident with  a vertex in $\Gamma i$ and a vertex in $\Gamma j$.

Let $\tilde G$ be a $\Gamma$-symmetric graph with quotient graph $G=(V,E)$. For each vertex orbit $\Gamma i$ we fix a representative vertex $i'\in \Gamma i$. 
For each edge orbit $\Gamma e$, choose a representative edge
and an ordering of its endpoints. Apply a group element so that the
first endpoint becomes the chosen representative $i'$ of its vertex
orbit. If $j'$ is the chosen representative of the vertex orbit of
the second endpoint, the resulting ordered edge has the form
$(i',\gamma j')$ for a unique $\gamma\in\Gamma$, since the action on
vertices is free. We orient the quotient edge $\Gamma e$ from
$\Gamma i$ to $\Gamma j$ and assign it the \emph{gain}
$\psi(\Gamma e)=\gamma$. For a loop, the chosen ordering of the
endpoints specifies its orientation.
This gain assignment $\psi: E \to \Gamma$ is well-defined and the pair $(G,\psi)$ is called the \emph{quotient $\Gamma$-labeled graph} of $\tilde G$. 
Conversely,  $\tilde G$ is called the \emph{covering graph} of $(G,\psi)$.

Note that the gain assignment $\psi$ depends on the choice of representative vertices. We may switch the gain assignment on the directed multigraph $G$ by choosing a different set of vertex orbit representatives. We regard two
gain assignments on the directed multigraph $G$ as equivalent if one can be obtained from the other by such a switching operation. The gain assignment $\psi$ also depends on the choice of orientation for each edge of $G$. If the orientation of an edge $\Gamma e$ in $G$ is reversed then the induced gain $\psi(\Gamma e)$ is replaced with $\psi(\Gamma e)^{-1}$.

In general, for a group $\Gamma$, a pair $(G,\psi)$ of a directed multigraph $G$ and  a map $\psi:E\to \Gamma$ is called a \emph{$\Gamma$-labeled graph} if it is the quotient $\Gamma$-labeled  graph of a $\Gamma$-symmetric graph $\tilde G$. Since $\tilde G$ is assumed to be simple, $(G,\psi)$ has no  parallel edges with the same gains or  loops with trivial gains. We regard two $\Gamma$-labeled graphs as equivalent if
they are derived from the same $\Gamma$-symmetric covering graph.

  We extend this notion to frameworks. 
  Let $\tau:\Gamma\to O(\mathbb R^d)$ be a homomorphism.
Since $\Gamma$ is finite, $\tau(\Gamma)$ is a finite subgroup of
$O(\mathbb R^d)$. For a $\Gamma$-symmetric graph $\tilde G$
with respect to $\theta$, a framework $(\tilde G,\tilde p)$ in
$\mathbb R^d$ is \emph{$\tau(\Gamma)$-symmetric} if
\begin{equation}
\label{eq:symmetric_func}
\tau(\gamma) \tilde p(i)=\tilde p(\theta(\gamma) i) \qquad \text{for all } \gamma\in \Gamma \text{ and all } i\in \tilde V.
\end{equation}
 
 Let $(\tilde G,\tilde p)$ be a $\tau(\Gamma)$-symmetric framework and let $G=(V,E)$ be the quotient $\Gamma$-labeled  graph of $\tilde G$. 
 Following the convention that $V$ is identified with the set $\{1,\dots, k\}$ of representative vertices, one can define the {\em quotient $\tau(\Gamma)$-labeled framework} as the triple $(G,\psi, p)$
with $p_i=\tilde{p}_i$ for each representative vertex $i$.
In general, a triple $(G,\psi, p)$ is called a {\em $\tau(\Gamma)$-labeled framework}.  
The {\em covering} of $(G,\psi,p)$  is the $\tau(\Gamma)$-symmetric framework $(\tilde G, \tilde p)$,
where $\tilde G$ is the covering graph of $(G,\psi)$ and $\tilde p$ is uniquely determined from $p$ by (\ref{eq:symmetric_func}).

Since $\Gamma$ is finite, we may choose orthogonal
coordinates in which all entries of the matrices $\tau(\gamma)$,
$\gamma\in\Gamma$, are algebraic over $\mathbb Q$.
We work in such coordinates throughout. This simplifies the definition
of genericity.
A $\tau(\Gamma)$-labeled framework $(G,\psi,p)$ over
$\mathbb R$ is \emph{generic} if its coordinates are algebraically
independent over $\mathbb Q$.
A $\tau(\Gamma)$-symmetric framework $(\tilde G,\tilde p)$ is called
\emph{$\tau(\Gamma)$-generic} if its quotient
$\tau(\Gamma)$-labeled framework $(G,\psi,p)$ is generic.

 \section[Symmetric local and global rigidity]{Symmetric local and global rigidity}\label{sec:symgr}

\subsection[Symmetry-forced local and global rigidity]{Forced-symmetric local and global rigidity}\label{subsec:forced_global}

We now introduce forced-symmetric local and global rigidity.

Let $\tilde G$ be a $\Gamma$-symmetric graph with respect to
$\theta$, and let $(\tilde G,\tilde p)$ be a $\tau(\Gamma)$-symmetric
framework in $\mathbb R^d$.  We say that $(\tilde G,\tilde p)$ is
\emph{$\tau(\Gamma)$-symmetric (locally) rigid} in $\mathbb R^d$ if
there exists an open neighbourhood $U$ of $\tilde p$ within the space
of configurations of $\tau(\Gamma)$-symmetric realisations of
$\tilde G$ such that, for every $\tilde q\in U$, equivalence of
$(\tilde G,\tilde q)$ and $(\tilde G,\tilde p)$ implies their congruence.
It is \emph{$\tau(\Gamma)$-symmetric globally rigid} in $\mathbb R^d$ if
every $\tau(\Gamma)$-symmetric framework $(\tilde G,\tilde q)$ in
$\mathbb R^d$ which is equivalent to $(\tilde G,\tilde p)$ is congruent
to $(\tilde G,\tilde p)$.

\begin{remark}\label{rem:equivalent_representations}
Our definition of symmetric global rigidity also covers the slightly
more general setting in which the representation may be replaced by
an equivalent one. Two representations
$\tau_1,\tau_2:\Gamma\to O(\mathbb R^d)$ are called \emph{equivalent}
if there exists an isometry $h$ of $\mathbb R^d$ such that
$h^{-1}\tau_1(\gamma)h=\tau_2(\gamma)$ for all $\gamma\in\Gamma$.
In this case we write $\tau_1\simeq\tau_2$. The graph action $\theta$
remains fixed. The following lemma shows that allowing equivalent
representations does not change the notion of forced-symmetric global rigidity.
\end{remark}

\begin{lemma}\label{lem:equivalent_representations}
Let $(\tilde G,\tilde p)$ be a $\tau(\Gamma)$-symmetric framework in
$\mathbb R^d$ with respect to $\theta$. Then $(\tilde G,\tilde p)$ is
$\tau(\Gamma)$-symmetric globally rigid if and only if every framework
$(\tilde G,\tilde q)$ in $\mathbb R^d$ which is
$\tau'(\Gamma)$-symmetric with respect to the same action $\theta$ for
some $\tau':\Gamma\to O(\mathbb R^d)$ with $\tau'\simeq\tau$, and
which is equivalent to $(\tilde G,\tilde p)$, is congruent to
$(\tilde G,\tilde p)$.
\end{lemma}
\begin{proof}
Suppose that $(\tilde G,\tilde p)$ is $\tau(\Gamma)$-symmetric globally
rigid, and let $(\tilde G,\tilde q)$ be an equivalent framework with
representation $\tau'\simeq\tau$. By definition, there exists an isometry $h$ such that
$\tau'(\gamma)=h^{-1}\tau(\gamma)h$ for all $\gamma\in\Gamma$.
Then, for every $v\in\tilde V$ and $\gamma\in\Gamma$,
\[
(h\circ\tilde q)(\gamma v)
=h\tau'(\gamma)\tilde q(v)
=\tau(\gamma)(h\circ\tilde q)(v).
\]
Thus $(\tilde G,h\circ\tilde q)$ is $\tau(\Gamma)$-symmetric.
Since $h$ is an isometry, this framework is equivalent to
$(\tilde G,\tilde p)$, and hence congruent to it by symmetric global
rigidity. Consequently, $(\tilde G,\tilde q)$ is also congruent to
$(\tilde G,\tilde p)$. 

The converse direction follows by taking $\tau'=\tau$.
\end{proof}

\begin{figure}[t]
\centering
\begin{minipage}[t]{0.32\linewidth}
\centering
\begin{tikzpicture}[very thick,scale=1]
\tikzstyle{every node}=[circle, draw=black, fill=white, inner sep=0pt, minimum width=5pt];
    \path (-0.5,1.5) node (p1) {} ;
    \path (0.5,1.5) node (p2)  {} ;
   \path (-1.049,-1.183) node (p11) {} ;
    \path (-1.55,-0.317) node (p22) {} ;
    \path (1.55,-0.317) node (p111) {} ;
    \path (1.049,-1.183) node (p222) {} ;
    
    \node [rectangle, draw=white, fill=white] (b) at (-0.8,1.75) {$\tilde p_1$};
     \node [rectangle, draw=white, fill=white] (b) at (0.8,1.75) {$\tilde p_6$};
    \node [rectangle, draw=white, fill=white] (b) at (-1.95,-0.35) {$\tilde p_2$};
      \node [rectangle, draw=white, fill=white] (b) at (1.95,-0.35) {$\tilde p_5$};
      \node [rectangle, draw=white, fill=white] (b) at (-1.4,-1.2) {$\tilde p_3$};
      \node [rectangle, draw=white, fill=white] (b) at (1.4,-1.2) {$\tilde p_4$};
    
      \draw (p1) -- (p2);
    \draw (p11) -- (p22);
    \draw (p111) -- (p222);
     \draw (p1) -- (p11);
    \draw (p11) -- (p111);
    \draw (p111) -- (p1);
     \draw (p2) -- (p22);
    \draw (p22) -- (p222);
    \draw (p222) -- (p2);
    
         \draw [dashed, thin] (0,-2) -- (0,2);
             \draw [dashed, thin] (30:2.5cm) -- (210:2.5cm);
             \draw [dashed, thin] (150:2.5cm) -- (330:2.5cm);
    \node [rectangle, draw=white, fill=white] (b) at (2.2,-1.5) {$\sigma$};         
       \node [rectangle, draw=white, fill=white] (b) at (-2.2,-1.5) {$\sigma'$};     
          \node [rectangle, draw=white, fill=white] (b) at (-0.4,-2) {$\sigma''$};  
             
\node [rectangle, draw=white, fill=white] (b) at (0,-2.5) {(a)};
        \end{tikzpicture}
\end{minipage}\hfill
\begin{minipage}[t]{0.32\linewidth}
\centering
\begin{tikzpicture}[very thick,scale=1]
\tikzstyle{every node}=[circle, draw=black, fill=white, inner sep=0pt, minimum width=5pt];
    \path (-0.5,-1.5) node (p1) {} ;
    \path (0.5,-1.5) node (p2)  {} ;
   \path (-1.049,1.183) node (p11) {} ;
    \path (-1.55,0.317) node (p22) {} ;
    \path (1.55,0.317) node (p111) {} ;
    \path (1.049,1.183) node (p222) {} ;

    \node [rectangle, draw=white, fill=white] (b) at (-0.8,-1.75) {$\tilde q_1$};
     \node [rectangle, draw=white, fill=white] (b) at (0.8,-1.75) {$\tilde q_6$};
    \node [rectangle, draw=white, fill=white] (b) at (-1.95,0.35) {$\tilde q_2$};
      \node [rectangle, draw=white, fill=white] (b) at (1.95,0.35) {$\tilde q_5$};
      \node [rectangle, draw=white, fill=white] (b) at (-1.4,1.2) {$\tilde q_3$};
      \node [rectangle, draw=white, fill=white] (b) at (1.4,1.2) {$\tilde q_4$};
    
      \draw (p1) -- (p2);
    \draw (p11) -- (p22);
    \draw (p111) -- (p222);
     \draw (p1) -- (p11);
    \draw (p11) -- (p111);
    \draw (p111) -- (p1);
     \draw (p2) -- (p22);
    \draw (p22) -- (p222);
    \draw (p222) -- (p2);
    
         \draw [dashed, thin] (0,-2) -- (0,2);
             \draw [dashed, thin] (30:2.5cm) -- (210:2.5cm);
             \draw [dashed, thin] (150:2.5cm) -- (330:2.5cm);
             
             \node [rectangle, draw=white, fill=white] (b) at (2.2,-1.5) {$\sigma'$};         
       \node [rectangle, draw=white, fill=white] (b) at (-2.2,-1.5) {$\sigma$};     
          \node [rectangle, draw=white, fill=white] (b) at (-0.4,-2) {$\sigma''$};  
             
\node [rectangle, draw=white, fill=white] (b) at (0,-2.5) {(b)};
        \end{tikzpicture}
\end{minipage}\hfill
\begin{minipage}[t]{0.32\linewidth}
\centering
\begin{tikzpicture}[very thick,scale=1]
\tikzstyle{every node}=[circle, draw=black, fill=white, inner sep=0pt, minimum width=5pt];
    \path (-0.5,-1.5) node (p1) {} ;
    \path (0.5,-1.5) node (p2)  {} ;
   \path (-1.049,1.183) node (p11) {} ;
    \path (-1.55,0.317) node (p22) {} ;
    \path (1.55,0.317) node (p111) {} ;
    \path (1.049,1.183) node (p222) {} ;

    \node [rectangle, draw=white, fill=white] (b) at (-0.8,-1.75) {$\tilde q_1$};
     \node [rectangle, draw=white, fill=white] (b) at (0.8,-1.75) {$\tilde q_6$};
    \node [rectangle, draw=white, fill=white] (b) at (-1.95,0.35) {$\tilde q_4$};
      \node [rectangle, draw=white, fill=white] (b) at (1.95,0.35) {$\tilde q_3$};
      \node [rectangle, draw=white, fill=white] (b) at (-1.4,1.2) {$\tilde q_5$};
      \node [rectangle, draw=white, fill=white] (b) at (1.4,1.2) {$\tilde q_2$};
    
      \draw (p1) -- (p2);
    \draw (p11) -- (p22);
    \draw (p111) -- (p222);
     \draw (p1) -- (p11);
    \draw (p11) -- (p111);
    \draw (p111) -- (p1);
     \draw (p2) -- (p22);
    \draw (p22) -- (p222);
    \draw (p222) -- (p2);
    
         \draw [dashed, thin] (0,-2) -- (0,2);
             \draw [dashed, thin] (30:2.5cm) -- (210:2.5cm);
             \draw [dashed, thin] (150:2.5cm) -- (330:2.5cm);
             
             \node [rectangle, draw=white, fill=white] (b) at (2.2,-1.5) {$\sigma$};         
       \node [rectangle, draw=white, fill=white] (b) at (-2.2,-1.5) {$\sigma'$};     
          \node [rectangle, draw=white, fill=white] (b) at (-0.4,-2) {$\sigma''$};  
             
\node [rectangle, draw=white, fill=white] (b) at (0,-2.5) {(c)};
        \end{tikzpicture}
\end{minipage}
\caption{Plane frameworks with dihedral symmetry of order $6$.
(b) is the reflection of (a) in the $x$-axis. (a) and (b) are congruent,
with representations $\tau$ and $\tau'\simeq\tau$, respectively.
(c) is $\tau(\Gamma)$-symmetric and equivalent, but not congruent, to (a).}
\label{fig:dihedralex}\label{fig:dihedralex0}
\end{figure}
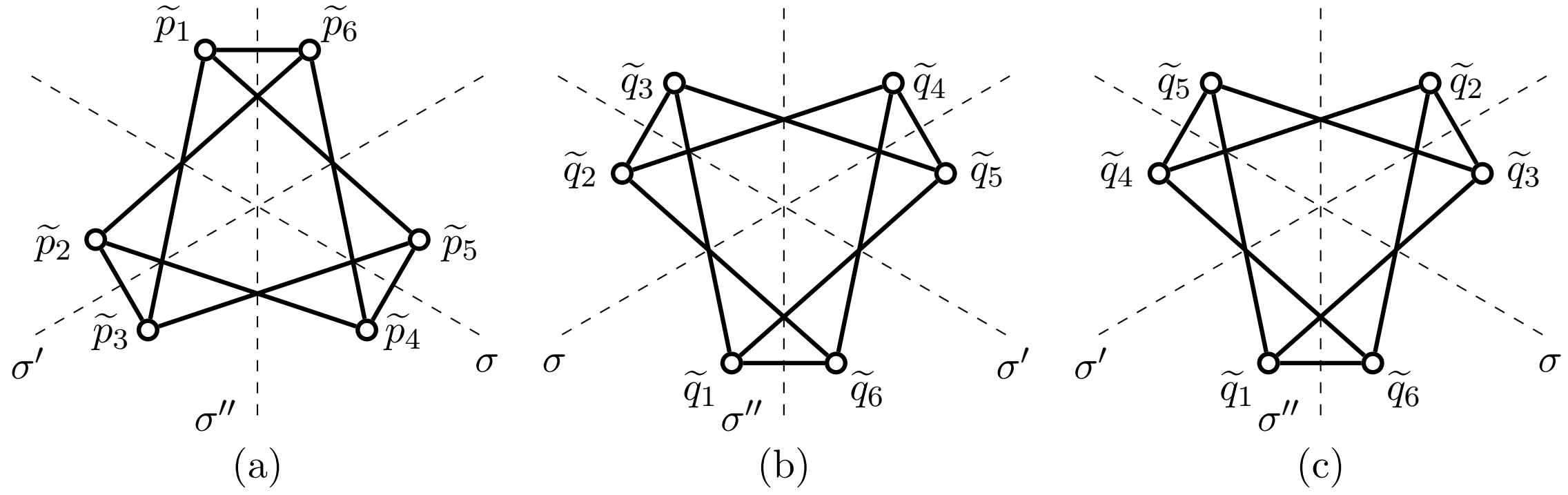

Consider, for example, the $\tau(\Gamma)$-symmetric plane framework $(\tilde G,\tilde p)$ in Figure~\ref{fig:dihedralex0}(a), where $\Gamma$ is the dihedral group of order 6. 
The frameworks in Figure~\ref{fig:dihedralex0}(b) and (c) are equivalent to $(\tilde G,\tilde p)$.
The framework in (c) is $\tau(\Gamma)$-symmetric  whereas
that in (b) is $\tau'(\Gamma)$-symmetric with 
$\tau'\simeq \tau$ since 
$h^{-1}\tau(\gamma)h=\tau'(\gamma)$ for all $\gamma\in\Gamma$,
where $h$ is the reflection in the $x$-axis.
The framework in (b) is congruent to $(\tilde G,\tilde p)$ since 
it can be obtained from $(\tilde G,\tilde p)$  by the reflection $h$.
On the other hand, the framework in (c) is not congruent to $(\tilde G,\tilde p)$
since the distance between $\tilde{p}_1$ and $\tilde{p}_2$
is different from that between $\tilde{q}_1$ and $\tilde{q}_2$.
Thus, the framework in (c) certifies that 
$(\tilde G,\tilde p)$ is not $\tau(\Gamma)$-symmetric globally rigid.

\subsection[Local and global rigidity of group-labeled graphs]{Local and global rigidity of group-labeled graphs}\label{subsec:labeled_global}

As explained in Section~\ref{sec:symbasic}, a
$\tau(\Gamma)$-symmetric framework can be represented compactly by
its quotient $\tau(\Gamma)$-labeled framework. It is therefore often
more convenient to analyse rigidity directly in terms of
$\tau(\Gamma)$-labeled frameworks. We now introduce notions of
local and global rigidity for $\tau(\Gamma)$-labeled frameworks that
correspond to $\tau(\Gamma)$-symmetric local and global rigidity,
and use this formulation
in what follows.

For  a $\Gamma$-labeled graph $(G,\psi)$ with $V=\{1,\ldots, k\}$ and a representation $\tau:\Gamma\to O(\mathbb R^d)$, we define the \emph{symmetric edge function}
$f_{G,\tau(\Gamma)}:\mathbb R^{dk}\rightarrow \mathbb R^{|E|}$ by
\[
f_{G,\tau(\Gamma)}(p)
=\bigl(\|p_i-\tau(\psi(e))p_j\|^2\bigr)_{e=ij\in E}
\qquad (p\in\mathbb R^{dk}).
\]

For a finite set $V$, define
$f_{V,\tau(\Gamma)}:\mathbb R^{d|V|}
\to\mathbb R^{|V|^2|\Gamma|}$ by
\[
f_{V,\tau(\Gamma)}(p)
=\bigl(\|p_i-\tau(\gamma)p_j\|^2\bigr)_{
(i,j,\gamma)\in V\times V\times\Gamma}.
\]

We say that a $\tau(\Gamma)$-labeled framework $(G,\psi,q)$ is \emph{equivalent} to another $\tau(\Gamma)$-labeled framework $(G,\psi,p)$ if $f_{G,\tau(\Gamma)}(p)=f_{G,\tau(\Gamma)}(q)$. Further, we say that $(G,\psi,p)$ and $(G,\psi,q)$ are \emph{congruent} if $f_{V,\tau(\Gamma)}(p)=f_{V,\tau(\Gamma)}(q)$.
This leads to the following notions of local and global
rigidity for $\tau(\Gamma)$-labeled frameworks.

A $\tau(\Gamma)$-labeled framework $(G,\psi,p)$ in
$\mathbb R^d$ is \emph{(locally) rigid} if there exists an open
neighbourhood $U$ of $p$ in $\mathbb R^{d|V|}$ such that, for every
$q\in U$, equivalence of $(G,\psi,q)$ and $(G,\psi,p)$ implies their
congruence. It is \emph{globally rigid} if every
$\tau(\Gamma)$-labeled framework $(G,\psi,q)$ in $\mathbb R^d$
which is equivalent to $(G,\psi,p)$ is congruent to $(G,\psi,p)$.

  \begin{lemma}
 \label{def:global2}
Let $(\tilde G,\tilde p)$ be a $\tau(\Gamma)$-symmetric
framework in $\mathbb R^d$ with respect to $\theta$, and let
$(G,\psi,p)$ be its quotient $\tau(\Gamma)$-labeled framework.
Then $(\tilde G,\tilde p)$ is $\tau(\Gamma)$-symmetric locally rigid
(respectively, globally rigid) if and only if $(G,\psi,p)$ is locally
rigid (respectively, globally rigid).
 \end{lemma}
  \begin{proof}
Consider a map $q\mapsto\tilde q$
that sends a configuration $q$ of $|V|$ points in $\mathbb R^d$
to a configuration $\tilde q$ of $|\tilde V|$ points in $\mathbb R^d$
by setting 
 \[
 \tilde q_{\gamma i}=\tau(\gamma)q_i \qquad (i\in V, \gamma\in\Gamma).
\]
This is a linear isomorphism from $\mathbb R^{d|V|}$ onto the space
of $\tau(\Gamma)$-symmetric configurations of $\tilde G$,
giving the configuration of the covering framework from a configuration of the quotient framework.
With this correspondence, the conditions for equivalence and congruence of the quotient frameworks are equivalent to the conditions for equivalence and congruence of the covering frameworks, respectively.
\end{proof}

\subsection[Symmetry-invariant isometries]{Symmetry-invariant isometries}\label{subsec:invariant_isometries}

For frameworks whose vertex orbit representatives affinely span
$\mathbb R^d$, local and global rigidity of the quotient framework can also
be described using isometries that preserve the prescribed symmetry.
An isometry $h$ of $\mathbb R^d$ is \emph{$\tau(\Gamma)$-invariant}
if $h\tau(\gamma)=\tau(\gamma)h$ for every $\gamma\in\Gamma$.
We denote the group of all such isometries by
$\mathrm{Euc}_{\tau(\Gamma)}$, the centraliser of $\tau(\Gamma)$
in $\mathrm{Euc}(d)$.

\begin{lemma}\label{lem:cs3}
Let $V$ be a finite set, let $p,q:V\to\mathbb R^d$ be two maps,
and let $\tau:\Gamma\to O(\mathbb R^d)$ be a homomorphism.
Suppose that $p(V)$ affinely spans $\mathbb R^d$.
Then $f_{V,\tau(\Gamma)}(p)=f_{V,\tau(\Gamma)}(q)$ if and only if
$q=h\circ p$ for some $h\in\mathrm{Euc}_{\tau(\Gamma)}$.
\end{lemma}
\begin{proof}
One direction is easy: if $q=h\circ p$ for a $\tau(\Gamma)$-invariant
isometry $h$, then for all $i,j\in V$ and $\gamma\in\Gamma$,
$
\|q_i-\tau(\gamma)q_j\|^2
=\|h(p_i)-h(\tau(\gamma)p_j)\|^2
=\|p_i-\tau(\gamma)p_j\|^2$.
Therefore $f_{V,\tau(\Gamma)}(p)=f_{V,\tau(\Gamma)}(q)$.

To see the converse, suppose that
$f_{V,\tau(\Gamma)}(p)=f_{V,\tau(\Gamma)}(q)$.
The coordinates corresponding to identity gains give
$\|p_i-p_j\|^2=\|q_i-q_j\|^2$ for all $i,j\in V$.
Thus there is an isometry $h$ such that $q_i=h(p_i)$ for every
$i\in V$.

We next show that $h$ is $\tau(\Gamma)$-invariant.
For all $i,j\in V$ and $\gamma\in\Gamma$,
\begin{equation}\label{eq:cs3}
\|p_i-\tau(\gamma)p_j\|^2
=\|q_i-\tau(\gamma)q_j\|^2
=\|p_i-h^{-1}(\tau(\gamma)h(p_j))\|^2.
\end{equation}
Fix $j$ and $\gamma$, and write $x=\tau(\gamma)p_j$ and
$y=h^{-1}(\tau(\gamma)h(p_j))$.
Expanding the squared distances in (\ref{eq:cs3}) and subtracting the equation for
any fixed $i_0\in V$ gives
$\langle p_i-p_{i_0},x-y\rangle=0$ for every $i\in V$.
Since $p(V)$ affinely spans $\mathbb R^d$, we obtain $x=y$.
Since this holds for all  $p_j$ and they affinely span $\mathbb R^d$, 
we obtain $h^{-1}\tau(\gamma)h=\tau(\gamma)$ for all $\gamma\in\Gamma$. 
\end{proof}

Let $(G,\psi,p)$ be a $\tau(\Gamma)$-labeled framework in
$\mathbb R^d$ with $p(V)$ affinely spanning $\mathbb R^d$.
Lemma~\ref{lem:cs3} implies that $(G,\psi,p)$ is locally rigid if and
only if there exists a neighbourhood $U$ of $p$ in
$\mathbb R^{d|V|}$ such that every equivalent
$\tau(\Gamma)$-labeled framework $(G,\psi,q)$ with $q\in U$
satisfies $q=h\circ p$ for some $h\in\mathrm{Euc}_{\tau(\Gamma)}$.
It is globally rigid if and only if every equivalent
$\tau(\Gamma)$-labeled framework $(G,\psi,q)$ in $\mathbb R^d$
satisfies $q=h\circ p$ for some $h\in\mathrm{Euc}_{\tau(\Gamma)}$.

\subsection[Characterising generic local rigidity]{Characterising generic local rigidity}\label{subsec:generic_local}

As in ordinary rigidity theory, generic local rigidity can be
characterised by the rank of the Jacobian of the symmetric edge
function, as follows. 

\begin{theorem}\label{thm:infequ}\label{thm:triv}
Let $(G,\psi,p)$ be a $\tau(\Gamma)$-generic
$\tau(\Gamma)$-labeled framework in $\mathbb R^d$ with $|V|\geq d$.
Then $(G,\psi,p)$ is locally rigid if and only if
\[
\operatorname{rank}df_{G,\tau(\Gamma)}(p)
=d|V|-\dim\mathrm{Euc}_{\tau(\Gamma)}.
\]
\end{theorem}
\begin{proof}
This statement follows by combining two known results through the notion of $\tau(\Gamma)$-symmetric infinitesimal rigidity.
Indeed, as in ordinary rigidity theory, $\tau(\Gamma)$-symmetric
local rigidity has a linearised counterpart, called
$\tau(\Gamma)$-symmetric infinitesimal rigidity, 
and local rigidity and infinitesimal rigidity 
coincide for $\tau(\Gamma)$-generic frameworks \cite{BSsiam}.
It is also known that the covering framework of $(G,\psi,p)$ is
$\tau(\Gamma)$-symmetric infinitesimally rigid if and only if the
Jacobian $df_{G,\tau(\Gamma)}(p)$ has rank $d|V|$ minus the dimension
of the space of trivial $\tau(\Gamma)$-symmetric infinitesimal
motions \cite{BSWWorbit}. Since $p$ is $\tau(\Gamma)$-generic and
$|V|\geq d$, this dimension equals
$\dim\mathrm{Euc}_{\tau(\Gamma)}$.
Together with Lemma~\ref{def:global2}, this proves the theorem.
\end{proof}

Since the Jacobian has the same rank at all $\tau(\Gamma)$-generic
realisations, Theorem~\ref{thm:infequ} implies that, for $|V|\geq d$,
local rigidity is a $\tau(\Gamma)$-generic property,
and it depends only on $\tau$ and the $\Gamma$-labeled graph $(G,\psi)$.

The dimension $\dim\mathrm{Euc}_{\tau(\Gamma)}$ depends only on
$\Gamma$ and its representation $\tau$. In dimensions 2 and 3,
it can be computed using standard character tables; see, for
example, \cite{Bis}.

\subsection[Global rigidity and connectivity]{Global rigidity and connectivity}\label{subsec:connectivity}

\begin{figure}[htp]\begin{center}
\begin{tabular}[b]{@{}c@{}}
\begin{tikzpicture}[very thick,scale=0.9]
\tikzstyle{every node}=[
    circle,
    draw=black,
    fill=white,
    inner sep=0pt,
    minimum width=5pt
];

\path (90:1.8cm)  node (p1) {};
\path (30:1.8cm)  node (p2) {};
\path (-30:1.8cm) node (p3) {};
\path (-90:1.8cm) node (p4) {};
\path (-150:1.8cm) node (p5) {};
\path (150:1.8cm) node (p6) {};

\node [rectangle, draw=white, fill=white]
    at (90:2.25cm) {$\tilde p_1$};
\node [rectangle, draw=white, fill=white]
    at (30:2.25cm) {$\tilde p_6$};
\node [rectangle, draw=white, fill=white]
    at (-30:2.25cm) {$\tilde p_5$};
\node [rectangle, draw=white, fill=white]
    at (-90:2.25cm) {$\tilde p_4$};
\node [rectangle, draw=white, fill=white]
    at (-150:2.25cm) {$\tilde p_3$};
\node [rectangle, draw=white, fill=white]
    at (150:2.25cm) {$\tilde p_2$};

\draw (p1) -- (p3);
\draw (p3) -- (p5);
\draw (p5) -- (p1);

\draw (p2) -- (p4);
\draw (p4) -- (p6);
\draw (p6) -- (p2);

\fill (0,0) circle (1pt);

       \end{tikzpicture}\\
(a)
\end{tabular}   
       \hspace{0.2cm}
\begin{tabular}[b]{@{}c@{}}
\begin{tikzpicture}[very thick,scale=0.9]
\tikzstyle{every node}=[
    circle,
    draw=black,
    fill=white,
    inner sep=0pt,
    minimum width=5pt
];

\path (-0.55, 2.0) node (p1) {};
\path ( 0.55, 2.0) node (p2) {};
\path (-0.55,-2.0) node (p3) {};
\path ( 0.55,-2.0) node (p4) {};

\path (-2.0, 0.55) node (p5) {};
\path (-2.0,-0.55) node (p6) {};
\path ( 2.0, 0.55) node (p7) {};
\path ( 2.0,-0.55) node (p8) {};

\node [rectangle, draw=white, fill=white]
    at (-0.85, 2.25) {$\tilde p_1$};
\node [rectangle, draw=white, fill=white]
    at ( 0.85, 2.25) {$\tilde p_2$};
\node [rectangle, draw=white, fill=white]
    at (-0.85,-2.25) {$\tilde p_3$};
\node [rectangle, draw=white, fill=white]
    at ( 0.85,-2.25) {$\tilde p_4$};

\node [rectangle, draw=white, fill=white]
    at (-2.35, 0.80) {$\tilde p_5$};
\node [rectangle, draw=white, fill=white]
    at (-2.35,-0.80) {$\tilde p_6$};
\node [rectangle, draw=white, fill=white]
    at ( 2.35, 0.80) {$\tilde p_7$};
\node [rectangle, draw=white, fill=white]
    at ( 2.35,-0.80) {$\tilde p_8$};

\draw (p1) -- (p2);
\draw (p2) -- (p4);
\draw (p4) -- (p3);
\draw (p3) -- (p1);

\draw (p1) -- (p4);
\draw (p2) -- (p3);

\draw (p5) -- (p7);
\draw (p7) -- (p8);
\draw (p8) -- (p6);
\draw (p6) -- (p5);

\draw (p5) -- (p8);
\draw (p6) -- (p7);

\draw[dashed,thin] (-2.7,0) -- (2.7,0);
\draw[dashed,thin] (0,-2.7) -- (0,2.7);

\draw[dashed,thin] (-2.3,-2.3) -- (2.3,2.3);
\draw[dashed,thin] (-2.3,2.3) -- (2.3,-2.3);

\node [rectangle, draw=white, fill=white]
    at (-2.45,-2.15) {$\sigma$};

\node [rectangle, draw=white, fill=white]
    at (2.45,-2.15) {$\sigma'$};

\end{tikzpicture}\\
(b)
\end{tabular}   
      \hspace{0.2cm}
\begin{tabular}[b]{@{}c@{}}
\begin{tikzpicture}[very thick,scale=0.9]
\tikzstyle{every node}=[
    circle,
    draw=black,
    fill=white,
    inner sep=0pt,
    minimum width=5pt
];

\path (-0.55, 2.0) node (p1) {};
\path ( 0.55, 2.0) node (p2) {};
\path (-0.55,-2.0) node (p3) {};
\path ( 0.55,-2.0) node (p4) {};

\path (-2.0, 0.55) node (p5) {};
\path (-2.0,-0.55) node (p6) {};
\path ( 2.0, 0.55) node (p7) {};
\path ( 2.0,-0.55) node (p8) {};

\node [rectangle, draw=white, fill=white]
    at (-0.85, 2.25) {$\tilde q_3$};
\node [rectangle, draw=white, fill=white]
    at ( 0.85, 2.25) {$\tilde q_4$};
\node [rectangle, draw=white, fill=white]
    at (-0.85,-2.25) {$\tilde q_1$};
\node [rectangle, draw=white, fill=white]
    at ( 0.85,-2.25) {$\tilde q_2$};

\node [rectangle, draw=white, fill=white]
    at (-2.35, 0.80) {$\tilde q_7$};
\node [rectangle, draw=white, fill=white]
    at (-2.35,-0.80) {$\tilde q_8$};
\node [rectangle, draw=white, fill=white]
    at ( 2.35, 0.80) {$\tilde q_5$};
\node [rectangle, draw=white, fill=white]
    at ( 2.35,-0.80) {$\tilde q_6$};

\draw (p1) -- (p2);
\draw (p2) -- (p4);
\draw (p4) -- (p3);
\draw (p3) -- (p1);

\draw (p1) -- (p4);
\draw (p2) -- (p3);

\draw (p5) -- (p7);
\draw (p7) -- (p8);
\draw (p8) -- (p6);
\draw (p6) -- (p5);

\draw (p5) -- (p8);
\draw (p6) -- (p7);

\draw[dashed,thin] (-2.7,0) -- (2.7,0);
\draw[dashed,thin] (0,-2.7) -- (0,2.7);

\draw[dashed,thin] (-2.3,-2.3) -- (2.3,2.3);
\draw[dashed,thin] (-2.3,2.3) -- (2.3,-2.3);
\node [rectangle, draw=white, fill=white]
    at (-2.45,-2.15) {$\sigma$};

\node [rectangle, draw=white, fill=white]
    at (2.45,-2.15) {$\sigma'$};

\end{tikzpicture}\\
(c)
\end{tabular} 
                 \end{center}
\vspace{-0.3cm}
\caption{Symmetric frameworks in the plane with disconnected graphs. 
Let $\mathcal{C}_n$ and $\mathcal{C}_{nv}$ denote 
the group of $n$-fold rotational symmetry and the dihedral group of order $2n$, respectively.  
(a) is $\mathcal{C}_6$-symmetric globally rigid. 
In this case, $\Gamma_C=\mathcal{C}_3$ for each connected component $C$.
On the other hand, 
(b) and (c) are not $\mathcal{C}_{4v}$-symmetric globally rigid
since they are equivalent but not congruent:
\(\|\tilde p_1-\tilde p_5\|\neq \|\tilde q_1-\tilde q_5\|\).
In this case, $\Gamma_C=\mathcal{C}_{2v}$ for each connected component $C$.
}
\label{fig:connectedex}
\end{figure}
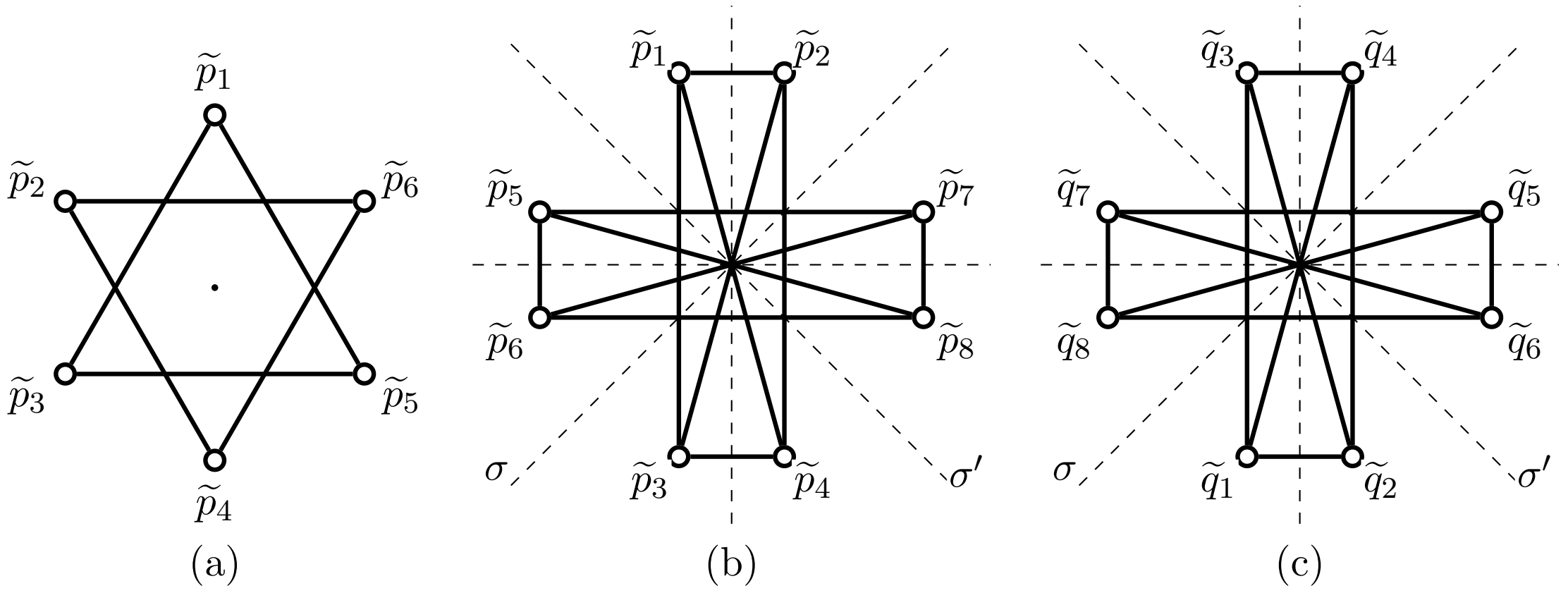

As pointed out in the introduction, 
a framework can be $\tau(\Gamma)$-symmetric globally rigid even if it is disconnected.
This is illustrated in Figure~\ref{fig:connectedex}.
In this subsection, we clarify the relationship between global rigidity and connectivity in more detail.

We first note that the connectivity of quotient $\tau(\Gamma)$-labeled graphs is necessary for global rigidity.

\begin{proposition}\label{prop:boundedness} Suppose that \((\widetilde G,\widetilde p)\) is a \(\tau(\Gamma)\)-generic framework which is \(\tau(\Gamma)\)-symmetric globally rigid in \(\mathbb R^d\). Then the quotient \(\Gamma\)-labeled graph \((G,\psi)\) is connected. \end{proposition}
\begin{proof}
 Suppose $G$ is not connected, and let \(X\) be a connected component of \(G\). Then we define a new quotient realisation \(q\) by \[ q_i= \begin{cases} -p_i,& i\in V(X),\\ p_i,& i\notin V(X). \end{cases} \] Since \(-I\) commutes with every element of \(\tau(\Gamma)\), all edge lengths within \(X\) are preserved, while all other edges are unchanged. Thus \(q\) gives a \(\tau(\Gamma)\)-symmetric realisation that is equivalent to \(p\). On the other hand, for \(i\in V(X)\) and \(j\notin V(X)\), \[ \|-p_i-p_j\|^2-\|p_i-p_j\|^2 = 4\langle p_i,p_j\rangle, \] which is non-zero by \(\tau(\Gamma)\)-genericity. Hence the corresponding covering frameworks are not congruent, contradicting \(\tau(\Gamma)\)-symmetric global rigidity. Thus \(G\) is connected.  
\end{proof}

Suppose that the quotient graph $G$ is connected, and let $C$ be a
connected component of $\tilde G$. The group $\Gamma$ acts transitively
on the connected components of $\tilde G$. Let
\[
\Gamma_C:=\operatorname{Stab}_{\Gamma}(C)=\{\gamma\in\Gamma:\gamma C=C\}.
\]
Equivalently, for any $v\in V(C)$, the subgroup $\Gamma_C$ consists of
the elements $\gamma$ for which $v$ and $\gamma v$ are connected by a
path in $\tilde G$. We may choose all vertex orbit representatives in
$C$. With this choice, every gain in $(G,\psi)$ 
belongs to $\Gamma_C$.

The next proposition describes how to reduce the problem of determining global rigidity of the entire framework to the problem of determining global rigidity of each connected component.

\begin{proposition}\label{prop:component_global}
Let $(\tilde G,\tilde p)$ be a $\tau(\Gamma)$-symmetric framework in
$\mathbb R^d$ whose quotient graph $G$ is connected. Let $C$ be a
connected component of $\tilde G$, and let
$\Gamma_C=\operatorname{Stab}_{\Gamma}(C)$. Suppose that
$\tilde p(V(C))$ affinely spans $\mathbb R^d$.
Then $(\tilde G,\tilde p)$ is $\tau(\Gamma)$-symmetric globally rigid
if and only if 
\begin{enumerate}
\item $(C,\tilde p|_{V(C)})$ is $\tau(\Gamma_C)$-symmetric globally rigid, and
\item $\mathrm{Euc}_{\tau(\Gamma_C)}=\mathrm{Euc}_{\tau(\Gamma)}$.
\end{enumerate}
\end{proposition}
\begin{proof}
Since the quotient graph $G$ is connected, we have
\begin{equation}\label{eq:cover}
  \{\delta v: \delta\in\Gamma, v\in V(C)\}=\tilde{V}.
\end{equation}

Suppose first that $(\tilde G,\tilde p)$ is
$\tau(\Gamma)$-symmetric globally rigid. 
To see that condition 1 holds, 
let $(C,r)$ be a $\tau(\Gamma_C)$-symmetric framework equivalent to $(C,\tilde p|_{V(C)})$.
We extend $r$ to a $\tau(\Gamma)$-symmetric framework $(\tilde G,\tilde q)$ by
\[
\tilde q(\delta v)=\tau(\delta)r(v)
\qquad (\delta\in\Gamma,\ v\in V(C)).
\]
We should verify that $\tilde q$ is well-defined.
Indeed, if $\delta v=\delta'v'$ for some $\delta,\delta'\in\Gamma$ and $v,v'\in V(C)$, 
then $\delta^{-1}\delta'C$ and $C$ share the vertex $v$ and hence
are the same connected component. Thus $\delta^{-1}\delta'\in\Gamma_C$.
Since $r$ is $\tau(\Gamma_C)$-symmetric,
$r(v)=r(\delta^{-1}\delta'v')=\tau(\delta^{-1}\delta')r(v')$,
and hence $\tilde q(\delta v)=\tau(\delta)r(v)=\tau(\delta')r(v')=\tilde q(\delta'v')$.

Thus, by (\ref{eq:cover}), $\tilde{q}$ is indeed 
a well-defined configuration of $\tilde{V}$.
Since $(C,r)$ is equivalent to $(C,\tilde p|_{V(C)})$, 
the extended framework $(\tilde G,\tilde q)$ is equivalent to $(\tilde G,\tilde p)$.
By global rigidity of $(\tilde G,\tilde p)$, 
$(\tilde G,\tilde p)$ is congruent to $(\tilde G,\tilde q)$.
The restriction to $C$ implies that $(C,\tilde p|_{V(C)})$ is congruent to $(C,r)$, proving condition 1.

To prove condition 2, let $h$ be a $\tau(\Gamma_C)$-invariant isometry
and extend $r=h\circ\tilde p|_{V(C)}$ to $\tilde q$ as above.
By global rigidity, $(\tilde G,\tilde q)$ is congruent to
$(\tilde G,\tilde p)$, so there is an isometry $g$ with
$\tilde q=g\circ\tilde p$ on all of $\tilde V$.
Since $g$ and $h$ agree on the affinely spanning set
$\tilde p(V(C))$, they are equal.
Thus, for every $\gamma\in\Gamma$ and $v\in V(C)$,
\[
h\tau(\gamma)\tilde p(v)
=g\tau(\gamma)\tilde p(v)
=g\tilde p(\gamma v)
=\tilde q(\gamma v)
=\tau(\gamma)\tilde q(v)
=\tau(\gamma)h\tilde p(v).
\]
Since $\tilde p(V(C))$ affinely spans $\mathbb R^d$, this implies
$h\tau(\gamma)=\tau(\gamma)h$, proving condition 2.

Conversely, suppose that conditions 1 and 2 hold for $(\tilde G,\tilde p)$, and let
$(\tilde G,\tilde q)$ be an equivalent $\tau(\Gamma)$-symmetric
framework. 
By condition 1, there is an isometry $h$ with $\tilde q(v)=h\tilde p(v)$ for every $v\in V(C)$.
For every $\gamma\in\Gamma_C$ and $v\in V(C)$, symmetry gives
\[
h\tau(\gamma)\tilde p(v)
=h\tilde p(\gamma v)
=\tilde q(\gamma v)
=\tau(\gamma)h\tilde p(v).
\]
Since $\tilde p(V(C))$ affinely spans $\mathbb R^d$, it follows that
$h\tau(\gamma)=\tau(\gamma)h$ for all $\gamma\in\Gamma_C$.
Thus $h$ is $\tau(\Gamma_C)$-invariant.
By condition 2, $h$ is $\tau(\Gamma)$-invariant, 
and hence for every $\delta\in\Gamma$ and $v\in V(C)$,
\[
\tilde q(\delta v)=\tau(\delta) \tilde{q}(v)=\tau(\delta)h\tilde p(v)
=h\tau(\delta)\tilde p(v)=h\tilde p(\delta v).
\]
By (\ref{eq:cover}) this implies $\tilde q(v)=h\tilde p(v)$ for all $v\in V(\tilde G)$.
\end{proof}

\section{Necessity of forced-symmetric redundant rigidity }\label{sec:Hen}

Suppose that   $(\tilde G,\tilde p)$ is $\tau(\Gamma)$-symmetric  rigid and the removal of any edge orbit still yields a $\tau(\Gamma)$-symmetric  rigid framework.
Then $(\tilde G,\tilde p)$ is called \emph{$\tau(\Gamma)$-symmetric redundantly rigid}. 
Equivalently,  the quotient $\tau(\Gamma)$-labeled framework $(G,\psi,p)$ is {\em  redundantly rigid} if the removal of any edge $e$ from $E(G)$ yields a  rigid $\tau(\Gamma)$-labeled framework.

The goal of this section is to prove the following forced-symmetric analogue of Hendrickson's theorem.
\begin{theorem}\label{thm:redundant_mirror} Let $(\tilde G,\tilde p)$ be a $\tau(\Gamma)$-generic framework in $\mathbb{R}^d$ 
 and  let $(G,\psi,p)$ be the quotient   $\tau(\Gamma)$-labeled framework. Suppose also that $|V(G)|\geq d+2$. If $(G,\psi,p)$ is  globally rigid, then it is redundantly rigid.
 \end{theorem}

Our proof strategy follows Hendrickson's original approach, 
and it uses the compactness of the configuration space under a suitable connectivity assumption.

\subsection{Extending Hendrickson's argument}\label{subsec:henext}

In order to prove Theorem~\ref{thm:redundant_mirror}, we consider
a more general setup of a polynomial map invariant under a group action.

For a field $K$ and an extension $K'$, let $td[K':K]$ denote the transcendence degree of the extension.
For a vector $p$ in $\mathbb{R}^k$, let $\mathbb{Q}(p)$ be the field generated by the entries of $p$ and the rationals.

Let \(X\) be a smooth manifold and \(f:X\rightarrow\mathbb{R}^m\) be a smooth map. 
Then \(x\in X\) is said to be a \emph{regular point} of \(f\) 
if the Jacobian \(Jf(x)\) has maximum rank, and is a \emph{critical point} of \(f\) otherwise. 
Also \(f(x)\) is said to be a \emph{regular value} of \(f\) 
if, for all \(y\in f^{-1}(f(x))\), \(y\) is a regular point of \(f\). 
Otherwise \(f(x)\) is called a \emph{critical value} of \(f\).

We consider the following setting.
Let $\mathcal{V}=\mathbb{R}^{d |V|}$ be the space of configurations of $|V|$ points in $\mathbb{R}^d$.
Let $\mathcal{V}_{\mathrm{full}}$ be the set of
full-dimensional configurations, that is, \[
\mathcal{V}_{\mathrm{full}}
=\{p\in\mathcal{V}:\operatorname{aff}p(V)=\mathbb R^d\}.
\]

Let $F$ be a closed subgroup of the $d$-dimensional Euclidean group $\mathrm{Euc}(d)$, 
and suppose that $F$ acts on $\mathcal{V}$ by the diagonal action.
Since the action of $F$ is free and proper over $\mathcal{V}_{\mathrm{full}}$, 
the quotient space $\mathcal{V}_{\mathrm{full}}/F$ is a smooth manifold.

Let $\mathbb{Q}'\subset \mathbb{R}$ be a finite extension of $\mathbb{Q}$.
Let $f:\mathcal{V}\rightarrow \mathbb{R}^m$ be a $\mathbb{Q}'$-polynomial map that is $F$-invariant.
Then, there is a well-defined smooth
map $\hat{f}:\mathcal{V}_{\mathrm{full}}/F\rightarrow \mathbb{R}^m$ such that $f|_{\mathcal{V}_{\mathrm{full}}}=\hat{f}\circ \pi$, 
where $\pi:\mathcal{V}_{\mathrm{full}}\rightarrow \mathcal{V}_{\mathrm{full}}/F$ is the quotient map.

We say that $p\in \mathcal{V}$ is \emph{$\mathbb{Q}'$-generic} if the entries of $p$ are algebraically independent over $\mathbb{Q}'$.

We will need the following lemmas. 
The first two are commonly used in rigidity theory, whereas the proof of the third one is technical and will be given in Appendix~\ref{subsec:lemma}.

\begin{lemma}\label{lem:edgefunction} (see e.g. \cite{jjsz}) 
Under the above setting, let $p$ be $\mathbb{Q}'$-generic.
Then
\[
td[\mathbb{Q}'(f(p)): \mathbb{Q}']\geq {\rm rank} Jf(p).
\]
\end{lemma}

\begin{lemma}\label{lem:regularvalue} (see e.g. \cite{gortler2010})
Under the above setting, let $p$ be $\mathbb{Q}'$-generic.
Then $f(p)$ is a regular value of $f$.
\end{lemma}

\begin{lemma}\label{lem:transcendence_degree_counting}
Under the above setting, let $p=(p_1,\dots, p_d)\in (\mathbb{R}^d)^d$. 
Then, there is $g\in F$ such that 
\[
td[\mathbb{Q}'(g\circ p):\mathbb{Q}']\leq d^2-\dim(F).
\]
\end{lemma}

Define $[m] = \{1,2,\ldots,m\}$.

\begin{theorem}\label{thm:hendricksongen}
Under the above setting, let $p$ be $\mathbb{Q}'$-generic and suppose further that 
\begin{enumerate}
\item $|V|\geq d+2$,
\item $\hat f^{-1}(\hat f([p]))=\{[p]\}$, and
\item for every $e\in [m]$, 
$f_{\bar{e}}^{-1}(f_{\bar{e}}(p))/F$ is compact,
where $f_{\bar{e}}$ denotes the map obtained from $f$ by removing the $e$-th coordinate.
\end{enumerate}
Then, for every $e\in [m]$, we have 
\[
{\rm rank} Jf_{\bar{e}}(p)=d|V|-\dim(F).
\]
\end{theorem}
\begin{proof}
By the assumption that $|V|\geq d+2$ and $p$ is $\mathbb{Q}'$-generic, 
$p\in \mathcal{V}_{\mathrm{full}}$ and hence $[p]\in \mathcal{V}_{\mathrm{full}}/F$.
By the second assumption, the fibre of $\hat f$ through
$[p]$ is a singleton. Since $p$ is $\mathbb Q'$-generic, we obtain
 \begin{equation}\label{eq:rank}
 \rank Jf(p)=d|V|-\dim(F).
 \end{equation}

Suppose for a contradiction that 
${\rm rank} Jf_{\bar{e}}(p)<d|V|-\dim(F)$ for some $e\in [m]$.
By \eqref{eq:rank}, and since deleting one row decreases the rank by at most
one, we have
\begin{equation}\label{eq:rank2}
\operatorname{rank} Jf_{\bar e}(p) = d|V| - \dim(F) - 1,
\end{equation}
and the $e$-th row of $Jf(p)$ is not in the span of the other rows.
Let $f_e$ denote the map obtained from $f$ by returning only the $e$-th
coordinate.

We claim the following:
\begin{claim}\label{claim:regular}
$f_{\bar{e}}^{-1}(f_{\bar{e}}(p))\subset \mathcal{V}_{\mathrm{full}}$.
\end{claim}
\begin{proof}
Suppose for a contradiction that there exists $q\in f_{\bar{e}}^{-1}(f_{\bar{e}}(p))$ such that $q\notin \mathcal{V}_{\mathrm{full}}$.
Then $q$ is contained in a proper affine subspace of $\mathbb{R}^d$.
Without loss of generality, assume that $\{q_1,\dots,q_d\}$ has the same affine span as  $\{q_1,\dots,q_{|V|}\}$.
By Lemma~\ref{lem:transcendence_degree_counting}, there exists $g\in F$ such that
\[
td[\mathbb{Q}'(g\circ (q_1,\dots,q_d)):\mathbb{Q}']\leq d^2-\dim(F).
\]
Since each of the points $q_i$ for $d+1\le i\le |V|$ is an affine combination of $q_1,\dots,q_d$, we have
\begin{equation}\label{eq:transcendence_degree1}
td[\mathbb{Q}'(g\circ q):\mathbb{Q}']\leq d^2-\dim(F)+(d-1)(|V|-d).
\end{equation}
On the other hand, by \eqref{eq:rank2} and Lemma~\ref{lem:edgefunction}, we have
\begin{equation}\label{eq:transcendence_degree2}
td[\mathbb{Q}'(f_{\bar{e}}(p)):\mathbb{Q}']\geq {\rm rank}\, Jf_{\bar{e}}(p)=d|V|-\dim(F)-1.
\end{equation}
Since $f_{\bar{e}}$ is $F$-invariant and $\mathbb{Q}'$-polynomial, we also obtain
\begin{equation}\label{eq:transcendence_degree3}
\begin{aligned}
td[\mathbb{Q}'(f_{\bar{e}}(p)):\mathbb{Q}']
&=td[\mathbb{Q}'(f_{\bar{e}}(q)):\mathbb{Q}']\\
&=td[\mathbb{Q}'(f_{\bar{e}}(g\circ q)):\mathbb{Q}']\\
&\leq td[\mathbb{Q}'(g\circ q):\mathbb{Q}'].
\end{aligned}
\end{equation}
Combining \eqref{eq:transcendence_degree1}, \eqref{eq:transcendence_degree2}, and \eqref{eq:transcendence_degree3}, 
we obtain $|V|\leq d+1$, contradicting the assumption that $|V|\geq d+2$.
\end{proof}

Claim~\ref{claim:regular} implies that $F$ acts freely and properly on $f_{\bar{e}}^{-1}(f_{\bar{e}}(p))$.
Moreover, by Lemma~\ref{lem:regularvalue}, $\hat{f}_{\bar{e}}([p])$ is a regular value of $\hat{f}_{\bar{e}}$.
Hence, (\ref{eq:rank2}) implies that  
$f_{\bar{e}}^{-1}(f_{\bar{e}}(p))/F$ is a 1-dimensional smooth embedded submanifold of $\mathcal{V}_{\mathrm{full}}/F$.
Since  $f_{\bar{e}}^{-1}(f_{\bar{e}}(p))/F$ is compact by the third assumption, 
it consists of a disjoint union of circles by the classification of 1-dimensional manifolds.
Let  \(\mathcal{O}\) be the component that contains $[p]$.

By (\ref{eq:rank2}), the function induced by $f_e$ on $\mathcal O$ has a non-zero differential at $[p]$,
and so there exist
  $[p']$ and $[p'']$ in \(\mathcal{O}\) with $f_{e}(p')< f_{e}(p)< f_{e}(p'')$.
  In particular, $[p']\neq [p'']$. 
  There are two arcs in $\mathcal O$ between $[p']$ and $[p'']$, one of which does not contain $[p]$, and so
 the Intermediate Value Theorem implies that there is a \([q]\in\mathcal{O}\) with \(f_{e}(q)=f_{e}(p)\)
 and \([q]\neq [p]\).
 \(f_{e}(p)=f_{e}(q)\) implies that $f(p)=f(q)$, contradicting the second assumption.
\end{proof}

\subsection{Symmetric version of Hendrickson's Theorem}\label{subsec:hend}

The symmetric version of Hendrickson's theorem given below is obtained
by applying Theorem~\ref{thm:hendricksongen} to the symmetric edge function defined in Section~\ref{subsec:labeled_global}.
In this setting, the compactness required by the third assumption of
Theorem~\ref{thm:hendricksongen} follows from the following lemma.

\begin{lemma}\label{lem:boundedness}
Suppose that $(\widetilde G,\widetilde p)$ is a $\tau(\Gamma)$-generic
framework which is $\tau(\Gamma)$-symmetric globally rigid in
$\mathbb R^d$. Suppose further that $d\geq 1$ and $|V|\geq d+2$.
Then, for every edge $e$ of the quotient
$\Gamma$-labeled graph $(G,\psi)$,
$
f^{-1}_{G-e}\bigl(f_{G-e}(p)\bigr)/\mathrm{Euc}_{\tau(\Gamma)}
$
is bounded.
\end{lemma}
\begin{proof}
Pick an edge \(e\) of \((G,\psi)\), and denote $\mathcal F_{\bar{e}}:=f_{G-e}^{-1}\bigl(f_{G-e}(p)\bigr)$.
We will show that \(\mathcal F_{\bar{e}}/\mathrm{Euc}_{\tau(\Gamma)}\) is bounded.

Let \(C\) be the vertex set of a connected component of \(\widetilde{G-e}\).
Recall that \(\Gamma_C=\operatorname{Stab}_{\Gamma}(C)\) denotes the
subgroup preserving \(C\), cf.~Section~\ref{subsec:connectivity}. 
We denote the translation subgroup of $\mathrm{Euc}_{\tau(\Gamma)}$ 
(resp. $\mathrm{Euc}_{\tau(\Gamma_C)}$) by $T_{\tau(\Gamma)}$ (resp. $T_{\tau(\Gamma_C)}$). 
Note that
\begin{equation}\label{eq:translation_subgroup}
T_{\tau(\Gamma)}=\{t\in \mathbb{R}^d:\tau(\gamma)t=t\text{ for all }\gamma\in\Gamma\}
\end{equation}
since $(t\circ \tau(\gamma))x =\tau(\gamma)x+t$ and 
$(\tau(\gamma)\circ t)x =\tau(\gamma)(x+t)$ for all $x\in \mathbb{R}^d$.

\begin{claim}\label{boundedmodF_C}
  There is a constant \(M_C>0\) such that, 
  for every \(q\in\mathcal F_{\bar{e}}\), there is a vector \(t_{q,C}\in T_{\tau(\Gamma_C)}\) such that
\(\|\tilde q_w-t_{q,C}\|\leq M_{C}\) for all \(w\in C\).
\end{claim}
\begin{proof}
The sum of the prescribed edge lengths in a spanning tree of \(C\) gives
a  bound \(M_C\) on the diameter of $\tilde{q}(C)$. 
Since each prescribed edge length is determined by $p$, 
$M_C$ is independent of the choice of $q\in\mathcal F_{\bar{e}}$.

Pick \(v\in C\) and let
\[
t_{q,C}:=\frac{1}{|\Gamma_C|}
\sum_{\gamma\in\Gamma_C}\tilde q_{\gamma v}
=\frac{1}{|\Gamma_C|}
\sum_{\gamma\in\Gamma_C}\tau(\gamma)\tilde q_v.
\]
This group average is fixed by every element of \(\tau(\Gamma_C)\), so
\(t_{q,C}\in T_{\tau(\Gamma_C)}\) by (\ref{eq:translation_subgroup}). 
Since \(\gamma v\in C\) for every \(\gamma\in\Gamma_C\),
we have, for every \(w\in C\),
\[
\|\tilde q_w-t_{q,C}\|
\leq\frac{1}{|\Gamma_C|}
\sum_{\gamma\in\Gamma_C}\|\tilde q_w-\tilde q_{\gamma v}\|
\leq M_C.
\]
\end{proof}

By Proposition~\ref{prop:boundedness}, \(G\) is connected, so \(G-e\) has
either one or two connected components. We treat these cases separately.

\smallskip
\noindent\emph{Case 1: \(G-e\) is connected:} 
Let \(C\) be the vertex set of a connected component of \(\widetilde{G-e}\).
By switching, we may assume that all vertex orbit representatives lie in \(C\).
Then, every gain of \(G-e\) belongs to \(\Gamma_C\). 

Our goal is to show that 
\begin{equation}\label{eq:boundedmodF}
T_{\tau(\Gamma_C)}=T_{\tau(\Gamma)}.
\end{equation}
Indeed, if (\ref{eq:boundedmodF}) holds, 
then $t_{q,C}$ shown in Claim~\ref{boundedmodF_C} 
belongs to $T_{\tau(\Gamma)}\subseteq \mathrm{Euc}_{\tau(\Gamma)}$ 
for each \(q\in\mathcal F_{\bar{e}}\).
So, translation by $t_{q,C}$ is a $\tau(\Gamma)$-invariant isometry, and 
the existence of the absolute bound \(M_C\) in Claim~\ref{boundedmodF_C} implies that
\(\mathcal F_{\bar{e}}/\mathrm{Euc}_{\tau(\Gamma)}\) is bounded.

To see  (\ref{eq:boundedmodF}),
suppose for a contradiction that \(T_{\tau(\Gamma_C)}\) is strictly larger than \(T_{\tau(\Gamma)}\).
Choose a vector \(t\in T_{\tau(\Gamma_C)}\setminus T_{\tau(\Gamma)}\).
By (\ref{eq:translation_subgroup}),  \(t\) can be chosen so that it has coordinates in the field $\mathbb{Q}(\tau(\Gamma))$. 
We define a configuration $q$ by
\[
q_i:=p_i+\varepsilon t\qquad(i\in V),
\]
where $\varepsilon\in \mathbb{R}$ is specified later.
Since every gain  of \(G-e\) is in $\Gamma_C$, 
(\ref{eq:translation_subgroup}) implies that 
$\|p_i-\tau(\psi(f))p_j\|=\|q_i-\tau(\psi(f))q_j\|$ for every edge \(f=ij\) of \(G-e\).
In other words, \(q\in\mathcal F_{\bar{e}}\) for any choice of \(\varepsilon\).

The next step is to choose \(\varepsilon\) so that 
the length of \(e=uv\) is also satisfied, 
i.e., \(\|q_u-\tau(\psi(e))q_v\|=\|p_u-\tau(\psi(e))p_v\|\).
For this, let  $a:=p_u-\tau(\psi(e))p_v$ and $b:=t-\tau(\psi(e))t$.
Then, the length constraint of $e$ is written as 
$\|a+\varepsilon b\|^2=\|a\|^2$.

If $b=0$, 
then $\|a+\varepsilon b\|^2=\|a\|^2$ is automatically satisfied for any choice of $\varepsilon$
and we choose any $\varepsilon\neq 0$.
If  $b\neq0$, then we choose $\varepsilon$ by
\[
\varepsilon=-\frac{2\langle a,b\rangle}{\|b\|^2},
\]
so that \(\|a+\varepsilon b\|^2=\|a\|^2\),
and $(G,\psi, q)$ is equivalent to $(G,\psi, p)$.
Since $(G,\psi,p)$ is $\tau(\Gamma)$-symmetric globally rigid, 
there is a \(\tau(\Gamma)\)-invariant Euclidean isometry \(h\) satisfying
$h(p_i)=q_i=p_i+\varepsilon t$ for all \(i\in V\).
Since \(p(V)\) affinely spans \(\mathbb R^d\), 
\(h\) must be the translation by \(\varepsilon t\), 
implying that \(\varepsilon t\in T_{\tau(\Gamma)}\). 

Since \(t\notin T_{\tau(\Gamma)}\), we must have \(\varepsilon=0\).
We have taken \(\varepsilon\neq 0\) if \(b=0\), so we must have \(b\neq 0\),
and hence 
\(\langle a,b\rangle =0\) by the definition of \(\varepsilon\). This contradicts the genericity of $p$, since $t$ has coordinates in
$\mathbb{Q}(\tau(\Gamma))$ and $b \neq 0$.
(If $e$ is a loop, note that
$\langle a, b\rangle = \langle p_u, w\rangle$ with
$w = (2I - \tau(\psi(e)) - \tau(\psi(e))^{-1})t$, and $w \neq 0$ since
$\langle t, w\rangle = \|b\|^2 \neq 0$.)
This  completes the proof of Case 1.

\smallskip
\noindent\emph{Case 2: \(G-e\) has two connected components.}
Write \(e=uv\), and let $V_1$ and $V_2$ be the vertex sets of the two
connected components of \(G-e\), with $u\in V_1$ and $v\in V_2$.
The covering of each connected component of $G-e$ 
consists of connected components of $\widetilde{G-e}$.
By switching, we may assume that $\psi(e)={\rm id}$ and all vertex orbit representatives in $V_i$ belong to the same connected component $C_i$ of $\widetilde{G-e}$ for each $i=1,2$.

Our goal is to show that
\begin{equation}\label{eq:boundedmodF2}
T_{\tau(\Gamma_{C_1})}=T_{\tau(\Gamma_{C_2})}=\{0\}.
\end{equation}
If this is the case, then by Claim~\ref{boundedmodF_C}
$\|\tilde q_i\|=\|q_i\|\leq \max\{M_{C_1}, M_{C_2}\}$ for all $i\in V$ and all $q\in\mathcal F_{\bar{e}}$.
Thus, \(\mathcal F_{\bar{e}}/\mathrm{Euc}_{\tau(\Gamma)}\) is bounded.

To see (\ref{eq:boundedmodF2}), suppose for a contradiction that \(T_{\tau(\Gamma_{C_1})}\) is non-zero.
Choose \(0\neq t\in T_{\tau(\Gamma_{C_1})}\)
with coordinates in the field $\mathbb{Q}(\tau(\Gamma))$.
We define a configuration $q$ by
\[
q_i:=
\begin{cases}
p_i+\varepsilon t, & (i\in V_1), \\
p_i, & (i\in V_2),
\end{cases}
\]
with 
\[
\varepsilon=-\frac{2\langle p_u-p_v,t\rangle}{\|t\|^2},
\]
which is non-zero since $p$ is generic.
Since the gain of each edge on $V_1$ belongs to $\Gamma_{C_1}$, 
(\ref{eq:translation_subgroup}) implies that 
$\|p_i-\tau(\psi(f))p_j\|=\|q_i-\tau(\psi(f))q_j\|$ for every edge \(f=ij\) in $G-e$.
Moreover, the choice of $\varepsilon$ guarantees that the length constraint of $e=uv$ is also satisfied.
Thus, $(G,\psi, q)$ is equivalent to $(G,\psi, p)$.

To derive a contradiction, we show that $(G,\psi, q)$ is not congruent to $(G,\psi, p)$.
If there is an \(x\in V_1\setminus\{u\}\), then
\begin{align*}
\|q_x-q_v\|^2-\|p_x-p_v\|^2
&=\|p_x+\varepsilon t-p_v\|^2-\|p_x-p_v\|^2 \\
&=2\langle \varepsilon t, p_x-p_v\rangle+\varepsilon^2\|t\|^2 \\
&=2\varepsilon \langle t, p_x-p_u\rangle\neq0
\end{align*}
by the genericity of $p$.
On the other hand, 
if \(V_1=\{u\}\), then \(|V|\geq d+2\) gives a vertex
\(y\in V_2\setminus\{v\}\), 
and the identical argument shows that 
\[
\|q_u-q_y\|^2-\|p_u-p_y\|^2
=2\varepsilon \langle t, p_v-p_y \rangle\neq0.
\]
This contradicts the global rigidity of $(G,\psi, p)$.
This contradiction completes the proof of Case 2, and hence the proof of the lemma.
\end{proof}

We are now ready to  prove the symmetric version of Hendrickson's theorem.

\begin{proof}[Proof of Theorem~\ref{thm:redundant_mirror}]
Let $f_{G,\tau(\Gamma)}$ be the symmetric edge function, and let
$f_{G-e}$ be obtained by deleting the coordinate corresponding to
$e\in E(G)$. The fibre
$f_{G-e}^{-1}(f_{G-e}(p))/\mathrm{Euc}_{\tau(\Gamma)}$
is bounded by Lemma~\ref{lem:boundedness}.
It is closed since $f_{G-e}$ is a polynomial map.
Hence, the fibre is compact.
By Theorem~\ref{thm:hendricksongen}, applied with
$F=\mathrm{Euc}_{\tau(\Gamma)}$,
\[
\operatorname{rank}df_{G-e}(p)
=d|V|-\dim\mathrm{Euc}_{\tau(\Gamma)}
\]
for every $e\in E$. By Theorem~\ref{thm:infequ}, $(G-e,\psi,p)$
is locally rigid for every $e\in E$, and therefore
$(G,\psi,p)$ is redundantly rigid.
\end{proof}

\subsection{Remarks on symmetric Hendrickson-type theorems and the generic property}\label{subsec:remarks}

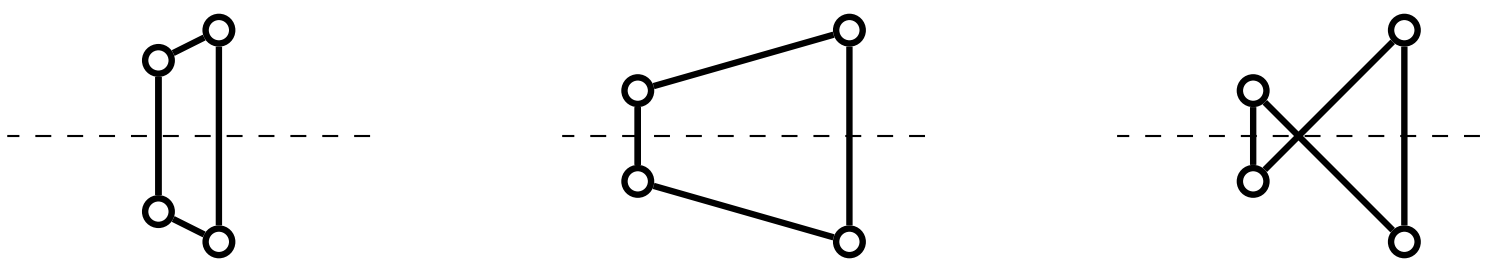
\begin{figure}[thp]
\begin{center}
         \begin{tikzpicture}[very thick,scale=1,rotate=90]
\tikzstyle{every node}=[circle, draw=black, fill=white, inner sep=0pt, minimum width=5pt];
    \path (-0.5,0.2) node (p1) {} ;
    \path (0.5,0.2) node (p2)  {} ;
    \path (-0.7,-0.2) node (p3)  {} ;
     \path (0.7,-0.2) node (p4)  {} ;
      \draw (p1) -- (p2);
    \draw (p3) -- (p4);
    \draw (p2) -- (p4);
    \draw (p1) -- (p3);
      \draw (p2) -- (p1);
         \draw [dashed, thin] (0,-1.2) -- (0,1.2);
         \node [rectangle, draw=white, fill=white] (b) at (-1.2,0) {(a)};
        \end{tikzpicture}
          \hspace{1cm}
         \begin{tikzpicture}[very thick,scale=1,rotate=90]
\tikzstyle{every node}=[circle, draw=black, fill=white, inner sep=0pt, minimum width=5pt];
    \path (-0.3,0.7) node (p1) {} ;
    \path (0.3,0.7) node (p2)  {} ;
    \path (-0.7,-0.7) node (p3)  {} ;
     \path (0.7,-0.7) node (p4)  {} ;
      \draw (p1) -- (p2);
    \draw (p3) -- (p4);
    \draw (p2) -- (p4);
    \draw (p1) -- (p3);
      \draw (p2) -- (p1);
         \draw [dashed, thin] (0,-1.2) -- (0,1.2);
         \node [rectangle, draw=white, fill=white] (b) at (-1.2,0) {(b)};
        \end{tikzpicture}
         \hspace{1cm}
         \begin{tikzpicture}[very thick,scale=1,rotate=90]
\tikzstyle{every node}=[circle, draw=black, fill=white, inner sep=0pt, minimum width=5pt];
    \path (-0.3,0.3) node (p1) {} ;
    \path (0.3,0.3) node (p2)  {} ;
    \path (-0.7,-0.7) node (p3)  {} ;
     \path (0.7,-0.7) node (p4)  {} ;
      \draw (p1) -- (p2);
    \draw (p3) -- (p4);
    \draw (p1) -- (p4);
      \draw (p2) -- (p3);
         \draw [dashed, thin] (0,-1.2) -- (0,1.2);
         \node [rectangle, draw=white, fill=white] (b) at (-1.2,0) {(c)};
        \end{tikzpicture}
                \end{center}
\vspace{-0.3cm}
\caption{(a), (b) Two $\mathcal{C}_s$-generic realisations of the cycle graph $C_4$ in $\mathbb{R}^2$. The framework in (a) is $\mathcal{C}_s$-symmetric globally rigid, whereas the framework in (b) is not, as illustrated in (c).}
\label{ex:Cs_example}
\end{figure}

The assumption $|V|\ge d+2$ in Theorem~\ref{thm:redundant_mirror}
is necessary when $\Gamma$ is trivial: a $d$-simplex has $d+1$
vertices and is globally rigid, but not redundantly rigid. 
For non-trivial symmetry groups, it remains open whether a smaller bound would suffice. 
Our current conjecture is that $|V|\geq d+1$ is sufficient for all non-trivial symmetry groups $\tau(\Gamma)$.
However, the exact bound may depend on the specific symmetry group $\tau(\Gamma)$.
For example, for the reflection group  in $\mathbb{R}^2$, $|V|\geq 3$ is necessary 
as shown in Figure~\ref{ex:Cs_example} (which will be explained in detail below),
whereas for the three-fold rotation group  in $\mathbb{R}^2$, $|V|\geq 2$ seems sufficient.
The proof of Theorem~\ref{thm:redundant_mirror} shows that the difficulty in lowering the bound arises 
when ensuring the free action of the group of
$\tau(\Gamma)$-invariant isometries on a generic fibre of the measurement map.

This seems to be also critical when extending the theorem of Gortler, Healy, and Thurston \cite{gortler2010} to the symmetric setting.
Connelly and Terrell initiated the use of group representation theory to block-diagonalise stress matrices of symmetric frameworks \cite{ConTer}. 
The symmetry-adapted block-decomposition of Oba and Tanigawa \cite{oba2023} makes the technique more formal.
The trivial-isotypic block of the decomposition would lead to a natural stress matrix counterpart for the forced-symmetric global rigidity.
It is therefore natural to try to adapt the stress-matrix theory of Connelly~\cite{connelly1982}, Gortler, Healy, and Thurston \cite{gortler2010} to this setting (as was done for periodic frameworks in \cite{dstt}).
However, a related issue arises in the symmetric setting when 
requiring the action of the group of
$\tau(\Gamma)$-invariant isometries to be free along a generic fibre of the forced-symmetric measurement map.

In fact, we have  a counterexample for the $\tau(\Gamma)$-generic property of $\tau(\Gamma)$-symmetric global rigidity.
Figure~\ref{ex:Cs_example} shows two  realisations of the same $\mathbb{Z}_2$-symmetric graph which are both generic with reflection symmetry, one of which is forced-symmetric globally rigid while the other is not.
Analogous examples can easily be constructed for higher dimensions.
We are, however, not able to find such an example in $\mathbb{R}^d$ which has more than $d$ vertex orbits, and we even propose the following conjecture.
\begin{conjecture}\label{conj}
For any non-trivial symmetry group $\tau(\Gamma)$ in $\mathbb{R}^d$, if the group acts freely on the vertex set, and the number of vertex orbits is at least $d+1$, 
then the $\tau(\Gamma)$-symmetric global rigidity  is a $\tau(\Gamma)$-generic property.
\end{conjecture}

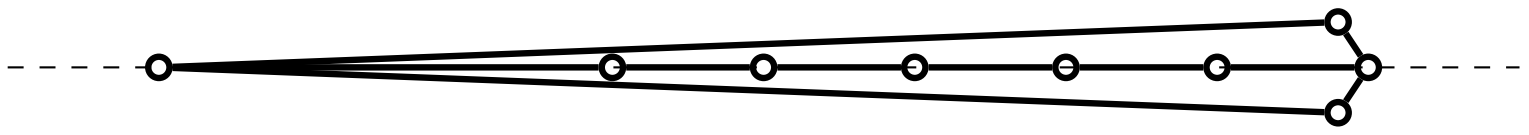
\begin{figure}[thp]
\begin{center}
         \begin{tikzpicture}[very thick,scale=1]
\tikzstyle{every node}=[circle, draw=black, fill=white, inner sep=0pt, minimum width=4pt];
    \path (0,0) node (p1) {} ;
    \path (-1,0) node (p2)  {} ;
    \path (-4,0) node (p5)  {} ;
     \path (1,0) node (q2)  {} ;
\path (2,0) node (q3) {} ;
    \path (3,0) node (q4)  {} ;
    \path (4,0) node (q5)  {} ;

\path (3.8,0.3) node (t)  {} ;
    \path (3.8,-0.3) node (b)  {} ;
     
      \draw (p1) -- (p2);
    \draw (p2) -- (p5);
         \draw (q2) -- (p1);
\draw (q2) -- (q3);
    \draw (q3) -- (q4);
      \draw (q4) -- (q5);
      \draw (t) -- (q5);
      \draw (t) -- (p5);
      \draw (b) -- (q5);
      \draw (b) -- (p5);
         \draw [dashed, thin] (-5,0) -- (5,0);
      
        \end{tikzpicture}
                        \end{center}
\vspace{-0.3cm}
\caption{A reflection-symmetric framework $(G,p)$ in $\mathbb{R}^2$ that is forced-symmetric globally rigid. However, there exists a reflection-symmetric framework $(G,q)$ in $\mathbb{R}^2$ which is not globally rigid, obtained by moving the top and bottom vertices farther from the mirror line, so that the vertices on the mirror line can be rearranged without changing edge lengths.}
\label{ex:Cs_example_non_free}
\end{figure}

Note that the above conjecture is not true in the case when the group action is not free on the vertex set, as shown in Figure~\ref{ex:Cs_example_non_free}.
Here, genericity is understood relative to the space of configurations satisfying the prescribed symmetry constraints.

\section{Body-bar frameworks} \label{sec:bb}

\subsection{Symmetric body-bar frameworks}

A body-bar framework consists of rigid bodies connected by bars.
We represent each body by a complete bar-joint framework containing
its attachment points.

Let $d\geq1$, and let $H$ be a loopless multigraph whose
vertices represent bodies and whose edges represent bars.
The associated \emph{body-bar graph} $G_H$ is the
simple graph constructed as follows.
\begin{itemize}
\item For each $v\in V(H)$, take $d+2$ core vertices
$x_{v,1},\ldots,x_{v,d+2}$. For each edge $e=uv\in E(H)$,
take two attachment vertices $x_{e,u},x_{e,v}$. All these vertices
are distinct.
\item The \emph{body} $B_H(v)$ is the complete graph on
\[
\{x_{v,1},\ldots,x_{v,d+2}\}
\cup\{x_{e,v}:e\in\delta_H(v)\},
\]
where $\delta_H(v)$ is the set of edges incident to $v$.
\item For each $e=uv\in E(H)$, add the bar
$x_{e,u}x_{e,v}$ between the corresponding bodies.
\end{itemize}
Thus $G_H$ consists of these complete bodies and
bars, and has $(d+2)|V(H)|+2|E(H)|$ vertices.
A pair $(G_H,p)$, where
$p:V(G_H)\to\mathbb R^d$, is called a
\emph{body-bar framework}.

To define the symmetric model, let $\Gamma$ be a finite group and
let $(H,\psi)$ be a $\Gamma$-labeled multigraph. We allow parallel
edges with equal gains and loops with non-trivial gains.
We first construct the associated quotient
\emph{$\Gamma$-labeled body-bar graph} $(G_H,\psi')$ as follows.
\begin{itemize}
\item For each $v\in V(H)$, take $d+2$ core vertices
$x_{v,1},\ldots,x_{v,d+2}$.
\item For each oriented edge $e=uv\in E(H)$, take two distinct
attachment vertices $a_e,b_e$, assigned to bodies $u,v$,
respectively. These are distinct vertices also when $u=v$.
All core and attachment vertices introduced in this construction
are distinct.
\item The \emph{body} $B_H(v)$ is the complete graph on the core
vertices of $v$ and all attachment vertices assigned to $v$.
Give every edge of $B_H(v)$ the identity gain.
\item For each $e=uv\in E(H)$, add the bar $a_e b_e$, directed
from $a_e$ to $b_e$, with gain $\psi(e)$.
\end{itemize}
Thus
$|V(G_H)|=(d+2)|V(H)|+2|E(H)|$.

\begin{figure}[htp]
\begin{center}
\begin{tikzpicture}[very thick,scale=1]
\tikzstyle{every node}=[circle, draw=black, fill=white, inner sep=0pt, minimum width=4pt];
      \node [draw=white, fill=white] (b) at (-0.7,0.75) {\small $\textrm{id}$};
     \node [draw=white, fill=white] (b) at (-0.2,0.75) {\small $\textrm{id}$};
\node [draw=white, fill=white] (b) at (0.7,0.75) {\small $\sigma$};
\node [draw=white, fill=white] (b) at (0.4,1.9) {\small $\sigma$};

   \path (0,0) node (u) [label =  below left: $u$] {};
\path (0,1.5) node (v) [label =  above left: $v$] {} ;
          \draw[->,>=stealth] (u) -- (v)  {} ;
        \path(v) edge [loop above,->, >=stealth,shorten >=2pt,looseness=45] (v);     
    \path (u) edge [->,>=stealth,bend right=60] (v);
     \path (u) edge [->,>=stealth,bend left=60] (v);

 \node [draw=white, fill=white] (b) at (0,-1) {(a)};
\end{tikzpicture}
\hspace{0.6cm}
\begin{tikzpicture}[very thick,scale=1]
\tikzstyle{every node}=[circle, draw=black, fill=white, inner sep=0pt, minimum width=4pt];

  \path (0,-0.2) node (u)  {};
  \path (-1,-0.5) node (uul){};
  \path (-0.1,-0.5) node (uur){};
  \path (-0.6,-0.7) node (uub){};
  \path (-0.8,0.2) node (uutl){};
  \path (-0.4,0.2) node (uutr){};

    \path (0.3,1.7) node (x)  {};
    \path (0.3,2.2) node (xx)  {};
  \path (-1,1.8) node (vvl){};
  \path (-0.8,2.3) node (vvr){};
  \path (-0.2,2.5) node (vvt){};
  \path (-0.8,1.2) node (vvbl){};
  \path (-0.4,1.2) node (vvbr){};
  
    \path (0.1,1.4) node (c){};

  \draw(uutl)--(vvbl);
  \draw(uutr)--(vvbr);

  \draw(u)--(uul);
  \draw(u)--(uur);
   \draw(u)--(uub);
  \draw(u)--(uutl);
    \draw(u)--(uutr);
  \draw(uul)--(uur);
  \draw(uul)--(uub);
  \draw(uul)--(uutl);
    \draw(uul)--(uutr);
   \draw(uur)--(uub);
  \draw(uur)--(uutl);
    \draw(uur)--(uutr);
    \draw(uub)--(uutl);
    \draw(uub)--(uutr);
     \draw(uutl)--(uutr);

  \draw(x)--(xx);
  \draw(x)--(vvl);
   \draw(x)--(vvr);
    \draw(x)--(vvt);
  \draw(x)--(vvbl);
   \draw(x)--(vvbr);
    \draw(x)--(c);
      \draw(xx)--(vvl);
   \draw(xx)--(vvr);
    \draw(xx)--(vvt);
  \draw(xx)--(vvbl);
   \draw(xx)--(vvbr);
    \draw(xx)--(c);
   \draw(vvl)--(vvr);
    \draw(vvl)--(vvt);
  \draw(vvl)--(vvbl);
   \draw(vvl)--(vvbr);
    \draw(vvl)--(c);
    \draw(vvr)--(vvt);
  \draw(vvr)--(vvbl);
   \draw(vvr)--(vvbr);
    \draw(vvr)--(c);
    \draw(vvbl)--(vvbr);
    \draw(vvbl)--(c);
      \draw(vvbr)--(c);
   \draw(vvbr)--(vvbl);
   \draw(vvt)--(xx);
  \draw(vvt)--(vvl);
   \draw(vvt)--(vvr);
    \draw(vvt)--(vvbl);
   \draw(vvt)--(vvbr);
    \draw(vvt)--(c);
 
  \path (u) edge[->,>=stealth] (c);
  \path (x) edge [->,>=stealth,bend right=80] (xx);
 \node [draw=white, fill=white] (b) at (0.3,0.7) {\small $\sigma$};
\node [draw=white, fill=white] (b) at (0.8,2.1) {\small $\sigma$};

 \node [draw=white, fill=white] (b) at (0,-1) {(b)};
\end{tikzpicture}
\hspace{0.6cm}
\begin{tikzpicture}[very thick,scale=1]
\tikzstyle{every node}=[circle, draw=black, fill=white, inner sep=0pt, minimum width=4pt];

 \node [draw=white, fill=white] (b) at (-1.6,2) {\small $B_H(v)$};
 \node [draw=white, fill=white] (b) at (-1.5,0) {\small $B_H(u)$};
\node [draw=white, fill=white] (b) at (3.6,2) {\small $\sigma B_H(v)$};
 \node [draw=white, fill=white] (b) at (3.5,0) {\small $\sigma B_H(u)$};

  \path (0,-0.2) node (u)  {};
  \path (-1,-0.5) node (uul){};
  \path (-0.1,-0.5) node (uur){};
  \path (-0.6,-0.7) node (uub){};
  \path (-0.8,0.2) node (uutl){};
  \path (-0.4,0.2) node (uutr){};

    \path (0.3,1.7) node (x)  {};
    \path (0.3,2.2) node (xx)  {};
  \path (-1,1.8) node (vvl){};
  \path (-0.8,2.3) node (vvr){};
  \path (-0.2,2.5) node (vvt){};
  \path (-0.8,1.2) node (vvbl){};
  \path (-0.4,1.2) node (vvbr){};
  
    \path (0.1,1.4) node (c){};

    \path (2,-0.2) node (u')  {};
  \path (2.1,-0.5) node (uul'){};
  \path (3,-0.5) node (uur'){};
  \path (2.6,-0.7) node (uub'){};
  \path (2.4,0.2) node (uutl'){};
  \path (2.8,0.2) node (uutr'){};

    \path (1.7,1.7) node (x')  {};
    \path (1.7,2.2) node (xx')  {};
  \path (3,1.8) node (vvl'){};
  \path (2.8,2.3) node (vvr'){};
  \path (2.2,2.5) node (vvt'){};
  \path (2.4,1.2) node (vvbl'){};
  \path (2.8,1.2) node (vvbr'){};
  
  \path (1.9,1.4) node (c'){};

  \draw(uutl)--(vvbl);
  \draw(uutr)--(vvbr);
  
  \draw(uutl')--(vvbl');
  \draw(uutr')--(vvbr');
 
\draw(u)--(c');
  \draw(u')--(c);
 
 \draw(x)--(xx');
  \draw(xx)--(x');

  \draw(u)--(uul);
  \draw(u)--(uur);
   \draw(u)--(uub);
  \draw(u)--(uutl);
    \draw(u)--(uutr);
  \draw(uul)--(uur);
  \draw(uul)--(uub);
  \draw(uul)--(uutl);
    \draw(uul)--(uutr);
   \draw(uur)--(uub);
  \draw(uur)--(uutl);
    \draw(uur)--(uutr);
    \draw(uub)--(uutl);
    \draw(uub)--(uutr);
     \draw(uutl)--(uutr);
     
      \draw(u')--(uul');
  \draw(u')--(uur');
   \draw(u')--(uub');
  \draw(u')--(uutl');
    \draw(u')--(uutr');
  \draw(uul')--(uur');
  \draw(uul')--(uub');
  \draw(uul')--(uutl');
    \draw(uul')--(uutr');
   \draw(uur')--(uub');
  \draw(uur')--(uutl');
    \draw(uur')--(uutr');
    \draw(uub')--(uutl');
    \draw(uub')--(uutr');
     \draw(uutl')--(uutr');

  \draw(x)--(xx);
  \draw(x)--(vvl);
   \draw(x)--(vvr);
    \draw(x)--(vvt);
  \draw(x)--(vvbl);
   \draw(x)--(vvbr);
    \draw(x)--(c);
      \draw(xx)--(vvl);
   \draw(xx)--(vvr);
    \draw(xx)--(vvt);
  \draw(xx)--(vvbl);
   \draw(xx)--(vvbr);
    \draw(xx)--(c);
   \draw(vvl)--(vvr);
    \draw(vvl)--(vvt);
  \draw(vvl)--(vvbl);
   \draw(vvl)--(vvbr);
    \draw(vvl)--(c);
    \draw(vvr)--(vvt);
  \draw(vvr)--(vvbl);
   \draw(vvr)--(vvbr);
    \draw(vvr)--(c);
    \draw(vvbl)--(vvbr);
    \draw(vvbl)--(c);
      \draw(vvbr)--(c);
   \draw(vvbr)--(vvbl);
   \draw(vvt)--(xx);
  \draw(vvt)--(vvl);
   \draw(vvt)--(vvr);
    \draw(vvt)--(vvbl);
   \draw(vvt)--(vvbr);
    \draw(vvt)--(c);
 
   \draw(x')--(xx');
  \draw(x')--(vvl');
   \draw(x')--(vvr');
    \draw(x')--(vvt');
  \draw(x')--(vvbl');
   \draw(x')--(vvbr');
    \draw(x')--(c');
      \draw(xx')--(vvl');
   \draw(xx')--(vvr');
    \draw(xx')--(vvt');
  \draw(xx')--(vvbl');
   \draw(xx')--(vvbr');
    \draw(xx')--(c');
   \draw(vvl')--(vvr');
    \draw(vvl')--(vvt');
  \draw(vvl')--(vvbl');
   \draw(vvl')--(vvbr');
    \draw(vvl')--(c');
    \draw(vvr')--(vvt');
  \draw(vvr')--(vvbl');
   \draw(vvr')--(vvbr');
    \draw(vvr')--(c');
    \draw(vvbl')--(vvbr');
    \draw(vvbl')--(c');
     \draw(vvbr')--(c');
   \draw(vvbr')--(vvbl');
   \draw(vvt')--(xx');
  \draw(vvt')--(vvl');
   \draw(vvt')--(vvr');
    \draw(vvt')--(vvbl');
   \draw(vvt')--(vvbr');
    \draw(vvt')--(c');

  \draw[dashed,thin](1,-0.8)--(1,2.5);
    \node [draw=white, fill=white] (b) at (0.7,-0.4) {\small $\sigma$};

 \node [draw=white, fill=white] (b) at (1,-1) {(c)};
\end{tikzpicture}
\end{center}
\vspace{-0.2cm}
\caption{(a) A $\mathcal C_s$-labeled multigraph
$(H,\psi)$, where $\mathcal C_s=\{\mathrm{id},\sigma\}$.
(b) The associated $\mathcal C_s$-labeled body-bar graph $(G_H,\psi')$.
(c) Its covering body-bar graph $\widetilde{G_H}$.
In (b) and (c), one core vertex per body and its incident edges are omitted.
}
\label{fig:bb}
\end{figure}
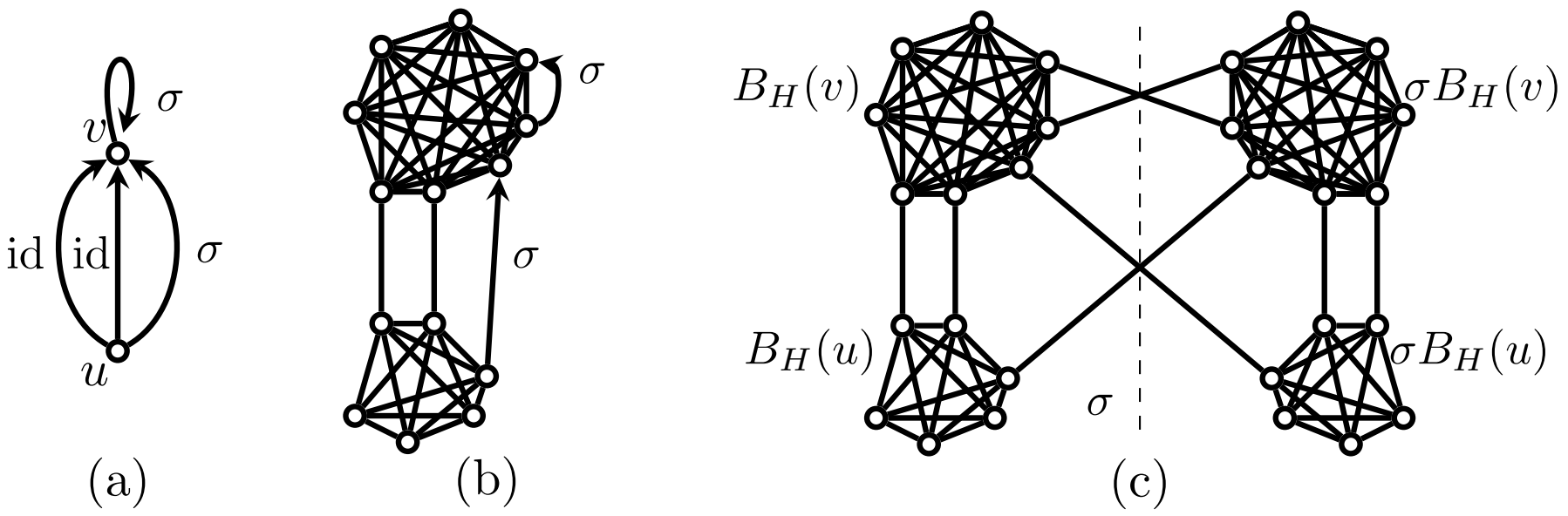

Let $\widetilde{G_H}$ be the covering graph of $(G_H,\psi')$.
Then, $\widetilde{G_H}$ is a body-bar graph in the ordinary sense.
See  Figure~\ref{fig:bb} for an example.

For a homomorphism $\tau:\Gamma\to O(\mathbb R^d)$ and a map
$p:V(G_H)\to\mathbb R^d$, the triple $(G_H,\psi',p)$ is called a
\emph{$\tau(\Gamma)$-labeled body-bar framework}. Its covering
framework $(\widetilde{G_H},\tilde p)$ is a
\emph{$\tau(\Gamma)$-symmetric body-bar framework}.

\subsection{Characterisation of generic forced-symmetric globally rigid body-bar frameworks}\label{sec:char_BB}

We now characterise generic forced-symmetric global rigidity of
body-bar frameworks in $\mathbb R^d$. 
We begin with the following auxiliary proposition.

\begin{proposition}\label{finiteprop}
Let $\Gamma$ be a finite group, let $(G,\psi)$ be a
$\Gamma$-labeled graph, and let
$\tau:\Gamma\to O(\mathbb R^d)$ be a homomorphism.
Let $K$ be the field generated by $\mathbb Q$ and the entries of
$\tau(\gamma)$ for $\gamma\in\Gamma$.
If $(G,\psi,p)$ is a $\tau(\Gamma)$-generic rigid framework in
$\mathbb R^d$, then there are only finitely many congruence classes
of $\tau(\Gamma)$-labeled frameworks equivalent to it.
Moreover, each of these classes contains a realisation whose
coordinates are algebraic over $K(p)$.
\end{proposition}
\begin{proof}
The maps $f_{G,\tau(\Gamma)}$ and $f_{V,\tau(\Gamma)}$ are
polynomial over $K$, and $f_{G,\tau(\Gamma)}$ is obtained by
selecting coordinates of $f_{V,\tau(\Gamma)}$.
The local rigidity of $(G,\psi,p)$ and the genericity of $p$ imply
$\operatorname{rank}df_{V,\tau(\Gamma)}(p)
=\operatorname{rank}df_{G,\tau(\Gamma)}(p)$.

Let $Y$ and $Z$ be the Zariski closures of the images of
$f_{G,\tau(\Gamma)}$ and $f_{V,\tau(\Gamma)}$, respectively.
The generic rank formula gives $\dim Z=\dim Y$, so the coordinate
projection $Z\to Y$ obtained by ignoring the coordinates not
corresponding to edges of $G$ is generically finite.
By the definition of congruence, distinct congruence classes of
realisations equivalent to $(G,\psi,p)$ give distinct points in the
fibre of the coordinate projection $Z\to Y$ over
$f_{G,\tau(\Gamma)}(p)$. Since $f_{G,\tau(\Gamma)}(p)$ is generic in
$Y$ over $K$, this fibre is zero-dimensional and hence finite.
Thus there are only finitely many such congruence classes.

The fibre is defined over $K(p)$, so its points have coordinates
algebraic over $K(p)$. Let $z$ be the point of this fibre
corresponding to one of these congruence classes.
A realisation $q:V\to\mathbb R^d$ belongs to this class precisely when
\[
f_{V,\tau(\Gamma)}(q)=z.
\]
This is a system of polynomial equations in the coordinates of $q$,
defined over $K(p,z)$, a finite algebraic extension of $K(p)$.
The system has a real solution, since the class is nonempty.
By the Tarski--Seidenberg transfer principle
\cite[Theorem~2.80]{bpr}, it has a real solution whose coordinates
are algebraic over $K(p,z)$, and hence over $K(p)$.
\end{proof}

\begin{theorem}\label{thm:BB}
Let $d\geq1$, let $\Gamma$ be a finite group, let $(H,\psi)$ be a
$\Gamma$-labeled multigraph with $V(H)\neq\varnothing$, and let
$\tau:\Gamma\to O(\mathbb R^d)$ be a homomorphism.
Then the following are equivalent:
\begin{enumerate}
\item Every $\tau(\Gamma)$-generic body-bar realisation of
$(G_H,\psi')$ is globally rigid in $\mathbb R^d$;
\item Some $\tau(\Gamma)$-generic body-bar realisation of
$(G_H,\psi')$ is globally rigid in $\mathbb R^d$;
\item $H$ is connected, $(G_H,\psi')$ is redundantly rigid in
$\mathbb R^d$, and
$
\mathrm{Euc}_{\tau(\Gamma_C)}
=\mathrm{Euc}_{\tau(\Gamma)}$,
where $C$ is a connected component of $\widetilde{G_H}$ and
$\Gamma_C=\operatorname{Stab}_{\Gamma}(C)$.
\end{enumerate}
\end{theorem}

\begin{proof}
The implication $1\Rightarrow2$ is immediate.

Suppose that condition 2 holds, and let $(G_H,\psi',p)$ be a
generic globally rigid realisation. Proposition~\ref{prop:boundedness}
implies that $H$ is connected. Each component $C$ of the covering graph contains a
full-dimensional body, so Proposition~\ref{prop:component_global}
gives $\mathrm{Euc}_{\tau(\Gamma_C)}
=\mathrm{Euc}_{\tau(\Gamma)}$.
Since
$|V(G_H)|=(d+2)|V(H)|+2|E(H)|\geq d+2$,
Theorem~\ref{thm:redundant_mirror} implies redundant rigidity.
Thus $2\Rightarrow3$.

To prove $3\Rightarrow1$, suppose now that condition 3 holds. Fix any generic realisation
$(G_H,\psi',p)$, and write $K=\mathbb Q(\tau(\Gamma))$ for the
field generated by the matrix entries of the representation.
For each realisation $q$ equivalent to $p$ and each $w\in V(H)$,
let $h_q^w$ be the unique isometry mapping $p$ to $q$ on $B_H(w)$.
This isometry exists and is unique because the body is complete
and its core vertices affinely span $\mathbb R^d$.

Let $e=uv$ be any edge of $H$.
The graph $G_{H-e}$ is obtained by deleting
the  two vertices $a_e, b_e$  from $G_H$. 
The redundant rigidity of $(G_H,\psi')$ implies that $(G_{H}-a_eb_e,\psi')$ is rigid.
Since each of $a_e$ and $b_e$ is adjacent to $d+1$ core vertices of its body, 
$(G_{H-e},\psi')$ is also rigid.

Let $p_e=p|_{V(G_{H-e})}$. Then $(G_{H-e},\psi',p_e)$
is a generic rigid realisation of $G_{H-e}$.
By Proposition~\ref{finiteprop}, we may choose a finite set
$\mathcal R(G_{H-e},p_e)$ of representatives of the realisations
equivalent to $(G_{H-e},\psi',p_e)$ modulo congruence, with
coordinates algebraic over $K(p_e)$.

For each $q_e\in\mathcal R(G_{H-e},p_e)$ and $w\in V(H)$, the
complete body $B_{H-e}(w)$ determines a unique isometry
\begin{equation}\label{eq:body_isometry}
h^w_{q_e}(x)=A^w_{q_e}x+t^w_{q_e}
\end{equation}
mapping $p_e$ to $q_e$ on that body. This uses the first $d+1$
core vertices, which are affinely independent and remain present
in $G_{H-e}$. 
Since $q_e$ is algebraic over $K(p_e)$, 
the entries of $A^w_{q_e}$ and $t^w_{q_e}$ are algebraic over $K(p_e)$.

By the genericity of $p$, the coordinates of $p(a_e),p(b_e)$
are algebraically independent over $K(p_e)$, and hence also over
the field obtained by adjoining the entries of all these isometries.

Now pick any realisation $q$ of $(G_H,\psi')$ equivalent to $p$.
Its restriction to $G_{H-e}$ is congruent to some
$q_e\in\mathcal R(G_{H-e},p_e)$.
By Lemma~\ref{lem:cs3}, after applying a $\tau(\Gamma)$-invariant
isometry to $q$, we may assume that this restriction equals $q_e$.
Then $h_q^w=h^w_{q_e}$ for each $w\in V(H)$, and in particular
\[
q(a_e)=A^u_{q_e}p(a_e)+t^u_{q_e},\qquad
q(b_e)=A^v_{q_e}p(b_e)+t^v_{q_e}.
\]
$G_H$ has an extra edge constraint for the edge $e=uv$,
which gives 
\[
\begin{aligned}
\|p(a_e)-\tau(\psi(e))p(b_e)\|^2
=\|A^u_{q_e}p(a_e)+t^u_{q_e}
-\tau(\psi(e))(A^v_{q_e}p(b_e)+t^v_{q_e})\|^2.
\end{aligned}
\]
By the algebraic independence of the two attachment points $p(a_e),p(b_e)$, this
is a polynomial identity in their coordinates. Comparing the
bilinear and linear terms gives
\[
\begin{aligned}
(A^u_{q_e})^T\tau(\psi(e))A^v_{q_e}&=\tau(\psi(e)),\\
(A^u_{q_e})^T(t^u_{q_e}-\tau(\psi(e))t^v_{q_e})&=0.
\end{aligned}
\]
Hence $A^u_{q_e}=\tau(\psi(e))A^v_{q_e}\tau(\psi(e))^{-1}$ and
$t^u_{q_e}=\tau(\psi(e))t^v_{q_e}$.
This implies 
\begin{equation}\label{eq:body_bar_isometries}
h^u_q=\tau(\psi(e))h^v_q\tau(\psi(e))^{-1}.
\end{equation}
This relation is preserved when we undo the normalisation, since
the applied isometry commutes with every element of $\tau(\Gamma)$.
Thus it holds for each $e\in E(H)$ and each equivalent realisation
$q$ of $(G_H,\psi')$.

Suppose, for the moment, that 
$H$ has no non-loop edges. Then $|V(H)|=1$, since $H$ is connected.
Let  $V(H)=\{u\}$.
Every equivalent realisation is obtained on $B_H(u)$ by a single
isometry $h^u_q$. For each loop $e$, (\ref{eq:body_bar_isometries}) gives
\begin{equation}\label{eq:loop_isometry}
  h^u_q\tau(\psi(e))=\tau(\psi(e))h^u_q.
\end{equation}
Since the gains of the loops $e$ generate $\Gamma_C$ for the component $C$ 
of the representative vertex, 
(\ref{eq:loop_isometry}) implies that 
$h^u_q$ commutes with every element of $\tau(\Gamma_C)$.
Thus, by the centraliser condition, 
$h^u_q$ commutes with every element of $\tau(\Gamma)$,
and hence $h^u_q$ is a $\tau(\Gamma)$-invariant isometry. Since $q=h^u_q\circ p$, it follows that $q$ is congruent to $p$.

Now suppose $H$  has at least one non-loop edge. Since $H$ is connected we can choose a spanning tree $T$. By switching we can assume that $\psi(e) = \mathrm{id}$ for every $e\in E(T)$. Now suppose that $e = uv \in E(T)$.
Then, by (\ref{eq:body_bar_isometries}), 
$h^u_q=h^v_q$ for every equivalent $q$, so the two bodies move
by the same isometry. 
Since this applies to every edge of $T$, all bodies move by a single isometry $h_q$. For every edge $e$ of $H$ that is not in $T$, (\ref{eq:body_bar_isometries}) then gives
$h_q\tau(\psi(e))=\tau(\psi(e))h_q$, and the gains of these edges generate
$\Gamma_C$. Hence, exactly as in the previous case, $h_q$ is a
$\tau(\Gamma)$-invariant isometry, and $q=h_q\circ p$ is congruent to $p$.
This completes the proof of $3\Rightarrow1$.
\end{proof}

\begin{remark}\label{rem:connected}
If the covering body-bar graph is connected, then $\Gamma_C=\Gamma$
and the centraliser condition is automatic. In this case,
Theorem~\ref{thm:BB} reduces to the characterisation by
redundant rigidity. 
\end{remark}

\begin{remark}
$\tau(\Gamma)$-symmetric  rigidity for $\tau(\Gamma)$-generic body-bar frameworks can be checked in polynomial time \cite{T2}. More precisely, since deleting an edge inside a body never destroys rigidity,
$(G_H,\psi')$ is redundantly rigid if and only if it is rigid and remains
rigid after deleting any bar $a_eb_e$, $e\in E(H)$, which requires at most
$|E(H)|+1$ rigidity tests. Moreover, after switching so that a spanning tree
of $H$ has identity gains, $\Gamma_C$ is generated by the gains of the
remaining edges, and for the fixed group $\Gamma$ and representation $\tau$
the centraliser condition depends only on the subgroup $\Gamma_C$.
Hence, for fixed $\Gamma$ and $\tau$,   Condition~3 of
Theorem~\ref{thm:BB} can be checked in  time polynomial in the size of $H$.
\end{remark}


\section{Further comments and future work}\label{sec:future}

There are several natural directions for future work. 
One  is to investigate forced-symmetric global rigidity for \emph{body-hinge} frameworks, consisting of full-dimensional rigid bodies in which each hinge connects exactly two bodies and allows one degree of freedom. The corresponding problem in the non-symmetric setting was solved in~\cite{jkt2}, and an extension to the forced-symmetric setting would provide a natural next step. Such results would also be relevant to the study of symmetric \emph{molecular} frameworks (i.e. body-hinge frameworks in which all hinge lines incident to each body are concurrent). A polynomial-time checkable combinatorial characterisation of generic forced-symmetric global rigidity of molecular frameworks, however, appears to be a more distant goal, since the corresponding  problem remains open even in the non-symmetric case.

Another natural question is whether forced-symmetric global rigidity is a generic property for bar-joint frameworks, provided we have at least $d+1$ vertex orbits (recall Conjecture~\ref{conj}). In particular, it would be interesting to obtain a combinatorial characterisation of forced-symmetric global rigidity
for $\tau(\Gamma)$-generic bar-joint frameworks in the plane. A classical Jackson--Jordán-type inductive approach \cite{jj},
combined with a symmetric version of the stress matrix (recall the discussion in Section~\ref{subsec:remarks}) or with the methods developed in~\cite{T1,kst}, appears promising,
although an additional difficulty arises from the need to handle $4$-regular circuits in the induction proof.

\section*{Acknowledgements}
We thank BIRS and CMI for supporting the BIRS-CMI Workshop on Rigidity, Flexibility and Complexity of Geometric Constraint Systems, where discussions related to this work took place.

Bernd Schulze was supported in part by UK Research and Innovation (grant number UKRI2397), under the EPSRC Mathematical Sciences Small Grant scheme. Csaba Kir\'aly was supported by the National Research,
Development and Innovation Office of Hungary, grant nos. PD-138102 and
Advanced 152786. 
Shin-ichi Tanigawa was supported by the Japan Science and Technology Agency (JST) as part of Adopting Sustainable Partnerships for Innovative Research Ecosystem (ASPIRE), Grant Number JPMJAP2520.

\paragraph{AI disclosure:} The bulk of this work was carried out before AI tools for proof discovery became available. Consequently, almost all results and proofs of this paper were obtained by the authors without the use of AI, with one exception. An initial proof of Lemma~\ref{lem:transcendence_degree_counting}  was first generated by ChatGPT. The  current proof given in Appendix~\ref{subsec:lemma} is a more elementary reformulation developed by the authors. Apart from that, AI was only used to check and polish some parts of the manuscript.

\appendix \section{Proof of Lemma~\ref{lem:transcendence_degree_counting}}\label{subsec:lemma}
In this section we provide a proof of Lemma~\ref{lem:transcendence_degree_counting}.

\begin{proof} 
For a configuration $q=(q_1,\ldots,q_d)$, define
\[
\Delta(q)
=
\begin{bmatrix}
q_2-q_1 & \cdots & q_d-q_1
\end{bmatrix}
\in M_{d,d-1}(\mathbb R).
\]
Put
\[
\begin{aligned}
s&=\operatorname{rank}\Delta(p)
=\dim\operatorname{aff}\{p_1,\ldots,p_d\}, \\
r&=\dim(F\cdot p).
\end{aligned}
\]

We first prove the lower bound (\ref{eq:dim}) of $r$.
Since $F$ is a closed subgroup of
$\mathrm{Euc}(d)$, its diagonal action is proper, and $F\cdot p$ is an
embedded submanifold diffeomorphic to $F/{\rm Stab}(F,p)$,
where 
\[
{\rm Stab}(F,p)
=
\{g\in F:gp_i=p_i\text{ for all }i=1,\ldots,d\}.
\]

Every element of ${\rm Stab}(F,p)$ fixes
$
{\rm aff}\{p_1,\ldots,p_d\}
$
pointwise,
and hence ${\rm Stab}(F,p)$ is a subgroup of the group of all Euclidean isometries fixing the
$s$-dimensional affine subspace ${\rm aff}\{p_1,\ldots,p_d\}$ pointwise, which is isomorphic to
$O(d-s)$. Hence
\[
\dim {\rm Stab}(F,p) \leq \dim O(d-s)=\binom{d-s}{2}.
\]
Therefore, we obtain
\begin{equation}\label{eq:dim}
r=\dim(F\cdot p)=\dim(F)-\dim({\rm Stab}(F,p))\geq \dim(F)-\binom{d-s}{2}.
\end{equation}

Next, consider the determinantal variety
\[
X_s
=
\left\{
q=(q_1,\ldots,q_d)\in(\mathbb R^d)^d:
\operatorname{rank}\Delta(q)\leq s
\right\}.
\]
It is an algebraic variety defined over $\mathbb Q$.
Since $q_1$ can be chosen freely and the
variety of $d\times(d-1)$ matrices of rank at most $s$ has dimension
$s(2d-1-s)$, 
\begin{equation}\label{eq:dimX}
\dim X_s
=d+s(2d-1-s)
=d^2-(d-s)(d-s-1).
\end{equation}

We claim that there is a linear map $
\ell:(\mathbb R^d)^d\longrightarrow\mathbb R^r
$
with rational coefficients such that
\begin{equation}\label{eq:iso}
\text{$d\ell\big|_{T_{p}(F\cdot p)}
:
T_{p}(F\cdot p)\longrightarrow\mathbb R^r$ is an isomorphism.}
\end{equation}
Indeed, since $\dim F\cdot p=r$, 
such a map exists with real coefficients,
and the full-rank condition is open, so its coefficients may be
chosen rational.

By the inverse function theorem, 
there is an open neighbourhood $U$ of $p$ in $F\cdot p$ 
that is mapped diffeomorphically by $\ell$ onto an open subset
$W\subseteq\mathbb R^r$. Choose $c\in W\cap\mathbb Q^r$, and let
$p'\in U$ be the point satisfying $\ell(p')=c$. Since
$p'\in F\cdot p$, there is $g\in F$ such that $p'=g\circ p$.

Define
\[
Z=X_s\cap\ell^{-1}(c).
\]
This is an algebraic set defined over $\mathbb Q$. 
Moreover, since $\ell\mid_U$ is a diffeomorphism and $\operatorname{rank}\Delta(p')=\operatorname{rank}\Delta(p)=s$ if $U$ is sufficiently small, 
$Z$ and $F\cdot p$ intersect transversally at $p'$ in $X_s$. 
In particular, $p'$ is a smooth point of $Z$ of local dimension $\dim X_s-r$.
Applying \cite[Lemma~2(c)]{JacksonMcCourtNixon2014} to a sufficiently small neighbourhood of $p'$, 
we obtain
\[
td[\mathbb{Q}(p'):\mathbb Q]\leq \dim X_s-r.
\]
Therefore, by (\ref{eq:dim}) and (\ref{eq:dimX}), we have
\[
\begin{aligned}
td[\mathbb Q(p'): \mathbb Q]
&\leq \dim X_s-r\\
&\mathrel{\textcolor{red}{\leq}}d^2-(d-s)(d-s-1)-\dim(F)+{d-s\choose 2}\\
&= d^2-\dim(F)-{d-s\choose 2}\\
&\leq d^2-\dim(F).
\end{aligned}
\]
Finally, since $\mathbb Q'/\mathbb Q$ is finite algebraic, adjoining
$\mathbb Q'$ does not change the transcendence degree. Thus
\[
td[\mathbb Q'(p'):\mathbb Q']
=td[\mathbb Q(p'): \mathbb Q]
\leq d^2-\dim(F),
\]
as required.
\end{proof}

\end{document}